\documentclass{amsart}
\usepackage{amsmath, amssymb, bm, formel,mathrsfs}
\newtheorem{thm}{Theorem}[section]
\newtheorem{lem}[thm]{Lemma}

\makeatletter

\@addtoreset{equation}{section}
\makeatother

\def\cal{\mathscr}

\def\bfa{{\mathbf a}}\def\bfc{{\mathbf c}}
\def\bfx{{\mathbf x}}\def\bfl{{\mathbf l}}
\def\bfX{{\mathbf X}}\def\bfm{{\mathbf m}}
\def\bfu{{\mathbf u}}
\def\bfy{{\mathbf y}}
\def\bfY{{\mathbf Y}}\def\bfW{{\mathbf W}}
\def\bfZ{{\mathbf Z}}
\def\bfz{{\mathbf z}}\def\bfj{{\mathbf j}}\def\bfs{{\mathbf s}}
\def\bfb{{\mathbf b}}\def\bfw{{\mathbf w}}
\def\bfu{{\mathbf u}}
\def\bfU{{\mathbf U}}
\def\bfv{{\mathbf v}}
\def\bfV{{\mathbf V}}
\def\bft{{\mathbf t}}\def\bfh{{\mathbf h}}
\def\Th{\Theta}
\def\bfone{{\mathbf 1}}
\def\om{\omega}
\def\dd{{\mathrm d}}
\def\MM{{\mathrm M}}
\begin{document}
\title[Counting in hyperbolic spikes]{Counting in hyperbolic spikes: 
the diophantine analysis of multihomogeneous
diagonal equations}
\begin{abstract}
A method  is described  to sum 
multi-dimensional arithmetic functions subject to hyperbolic summation
conditions, provided that asymptotic formulae in rectangular boxes are
available. In combination with the circle method, the new method is a
versatile tool to count rational points on
algebraic varieties defined by multi-homogeneous diagonal
equations.
\end{abstract}

\author{ Valentin Blomer and J\"org Br\"udern}
\address{Mathematisches Institut, Bunsenstrasse 3--5, D-37073 G\"ottingen, Germany}
\email{blomer@uni-math.gwdg.de, bruedern@uni-math.gwdg.de}

\thanks{First author supported by the Volkswagen Foundation and a Starting Grant
of the European Research Council.}

\maketitle
\section{Introduction}

Our main concern in this memoir is with the distribution of rational points on
a class of multiprojective varieties that we now introduce. Let $d,k,n$ be natural numbers with 
$n\ge 2$. Then, whenever $a_0,a_1,\ldots,a_n$ are non-zero integers, the equation
\be{01} \sum_{j=0}^n a_j (x_{1,j}x_{2,j} \cdots x_{k,j})^d = 0
\ee
is homogeneous of degree $d$ in the variables $\bfx_i = (x_{i,0},\ldots,x_{i,n})$,
and therefore defines a variety $\cal V$ in $\mathbb{P}^n(\QQ)^k$. Its rational points 
are in $1$-to-$2^k$ correspondence to solutions of \rf{01} in primitive vectors
$\bfx_i\in\ZZ^{n+1}$, with each $\bfx_i$ unique up to sign. Since the norm
$$ |\bfx_i|= \max_{0\le j\le n} |x_{i,j}| $$
is independent of the ambiguous sign, we may define the height of the rational point 
as
\be{02}  (|\bfx_1||\bfx_2|\cdots |\bfx_k|)^{n+1-d}. \ee
Motivated by an influential set of conjectures put forward by Manin (see \cite{FMT}), Batyrev and Tschinkel \cite{BT} and Peyre \cite{P}, we seek to establish an asymptotic formula for the number $\mathrm{N}(B)$ of rational points on $\cal V$ with $x_{i,j}\neq 0$ for all $1\le i\le k$, $0\le j\le n$ and height not exceeding $B$. Our methods prove successful whenever $n$ is sufficiently large, the threshold being determined by our understanding of moments of classical Weyl sums. In this context, let $n_0(d)$ denote the smallest even natural number with the property that for any positive real number $\eps$ one has
\be{03} \int_0^1 \Big| \sum_{1\le x\le P} {\mathrm e}^{2\pi{\mathrm i}\al x^d}\Big|^{n_0(d)}\,\dd\al \ll P^{n_0(d)-d +\eps}. \ee
The integral here has an interpretation as the number of solutions of an underlying diophantine equation, and in this way one finds that the integral is bounded below by $P^{n_{0}/2}$. This implies that $n_0(d)\ge 2d$. It is also immediate that $n_0(1)=2$ and $n_0(2)=4$, and as 
an averaged version of Hardy and Littlewood's well-known conjecture K, it is expected that $n_0(d)=2d$ holds for all $d$. The current records are
$$ n_0(3)\le 8, \quad n_0(4)\le16, \quad n_0(5)\le 28,\quad n_0(6)\le 44, $$
see \cite[Lemma 2.5]{hlm} for $d=3$ or $4$, and \cite{Wmulti} for $d= 5$ and $6$. When $d$ is large,
one has  $n_0(d) \le 1.543\,d^2$ (see \cite{Wmulti}).

\begin{thm}\label{Thm01} Let $d,k,n$ be natural numbers with $n\ge n_0(d)$. Let $a_0,\ldots,a_n$ be non-zero integers. Then, there are a positive number $\delta$, a non-negative number $C$, and a monic polynomial $Q\in\RR[X]$ of degree $k-1$ such that
$$\mathrm{N} (B) = CBQ(\log B) + O(B^{1-\delta}). $$
Moreover, the number $C$ is positive if and only if the equation
\be{04} a_0y_0^d + a_1y_1^d + \ldots + a_n y_n^d =0 \ee
has non-trivial solutions in real numbers and in $p$-adic numbers, for all primes $p$.
\end{thm}

The constant $C$ is a product of local densities and coincides with the predictions stemming from a formal use of the Hardy-Littlewood method. We shall make this more precise at the very end of this paper.

The counting function $\mathrm{N}(B)$ is blind for points on coordinate hyperplanes because the intersection
of the union of them with $\cal V$ contains abnormally many points of small height. To see this, first
suppose that $k\ge 3$, choose
$$ \bfx_1 =(1,0,0,\ldots,0), \quad \bfx_2 = (0,1,0,\ldots,0) $$
and $\bfx_i=\bfx_2$ for $3\le i\le k-1$. Then any primitive $\bfx_k\in\ZZ^{n+1}$ yields a point 
$(\bfx_1,\ldots,\bfx_k)$ on $\cal V$, and there are $\gg B^{(n+1)/(n+1-d)}$ such points where the height
does not exceed $B$. Similarly, when $k=2$, choose $(z_1,\ldots,z_n)\in\ZZ^n$ primitive and
$$ \bfx_1 =(1,0,0,\ldots,0), \quad \bfx_2 = (0,z_1,\ldots,z_n)$$
to find $B^{n/(n+1-d)}$ points on $\cal V$ with height not exceeding $B$. When $d>1$, the order of magnitude
here is still bigger than the leading term in Theorem \ref{Thm01}.

Only for very few choices for the parameters $d,k,n$ the conclusions in Theorem \ref{Thm01} are already known. 
Of course
when $k=1$, the equation \rf{01} becomes a standard diagonal equation, and the height \rf{02} reduces to
the naive height. In this case, an asymptotic  evaluation of $\mathrm{N}(B)$ is possible by classical versions of the circle method,
and there is an established theory that has been developed in parallel with
the asymptotic formula in Waring's problem (see Vaughan \cite{hlm} and Wooley \cite{W,Wmulti} for an account), 
and to which we have nothing to add. In contrast, all cases where $k\ge 3$ seem to be the first results 
at all about this family of varieties. In fact, there are remarkably few examples of multiprojective varieties where the rational points have been counted satisfactorily with respect to the anticanonical height
(see the authors' note \cite{BBPaddy}, Le Boudec \cite{Bou}, and Schindler \cite{Sch}  for recent efforts in a biprojective setting).  
For the varieties under consideration,
this height is given  by \rf{02}. It is the geometrically most natural and the analytically most challenging choice. Indeed, if one considers \rf{01} as an affine equation and cuts out a portion of its integral solutions by a constraint like $|\bfx_1||\bfx_2|\ldots|\bfx_k| \le X$, then the underlying surface has hyperbolic spikes where standard counting routines tend to deny service. 
In the special case $k=2$, $d=1$, a number of devices have been developed  
to overcome this difficulty. 
Bump \cite[Chapter 5]{Bump} and Vinogradov and Takhtadzhyan \cite{VT} independently observed that there is a
natural bijection between the rational points on
\be{05} 
x_0y_0 + x_1y_1 + x_2y_2 =0 
\ee
and the cosets of $\mathrm{SL}_3(\ZZ)$ factored by the group of upper triangular unipotent matrices. With this in hand,
the Dirichlet series
\be{06}
\sum_{\bfx,\bfy} (x_0^2+x_1^2+x_2^2)^{-s}(y_0^2+y_1^2+y_2^2)^{-s},
\ee
with the sum running over primitive solutions $\bfx\in\ZZ^3$, $\bfy\in\ZZ^3$ of \rf{05}, can be expressed in terms of the minimal parabolic Eisenstein series for $\mathrm{GL}_3$, and one finds
that the analytic function defined by \rf{06} extends to a meromorphic function on the whole complex plane.
In $\Re s >\f12$ its only poles are at $s=1$ and $s=\f34$ (see \cite[Theorem 7.1]{Bump}). Some mundane 
analytic number theory then yields a version of Theorem \ref{Thm01} for the equation \rf{05}, but with the euclidean
norm used in the definition of height.  This approach rests on the observation that the biprojective variety defined by \rf{05} carries a natural group structure. A similar line of thought is present in much work related to the Manin-Peyre conjecture, 
following the pivotal analysis of flag varieties by Franke, Manin and Tschinkel \cite{FMT}. In the higher degree cases of Theorem \ref{01}, however, an analogous theory is not available, and we will have to turn to other methods. Thunder \cite{Th} recovered the results of  Franke, Manin and Tschinkel by using ideas from the geometry of numbers, but again it appears difficult to apply  his method to the higher degree cases
in Theorem \ref{Thm01}. For another method to count primitive solutions of \rf{05} see Browning \cite{TB}.

It seems natural to appoach the equation \rf{01} directly by the circle method. A first attempt was made by Robbiani.
He also studied the bilinear case $d=1$, $k=2$ and used a version of the circle method that relies on an identity of Duke, Friedlander and Iwaniec \cite{dfi}, and that Heath-Brown \cite{HB,HB2} applied to the analytic theory of quadratic and cubic forms. Robbiani's result only covers the cases $n\ge 3$, and does not easily generalise to
higher degree. A more classical approach has been engineered by Spencer \cite{S}. His method is in the spirit of Vinogradov's description of the circle method. Spencer concentrates on the equation \rf{05} where the height contraint 
$|x_iy_j|^2\le B$ $(0\le i,j\le 2)$ prevents the  generating function
from immediate factorisation.
This causes considerable complication in detail. As Spencer shows, it is possible to disentangle the height condition with an appropriate Fourier analysis. However,  certain difficulties remain, and if $\mathrm{N}^*(B)$ denotes the function $\mathrm{N}(B)$
in the special case \rf{05}, then Spencer only shows that $\mathrm{N}^*(B) = C B\log B +O(B)$, with some unspecified constant $C$.

In this memoir we propose another approach to problems in this class that is based on a straightforward use of the circle method on the one hand, and on the other on a version of Dirichlet's method of the hyperbola in  weighted setting. Once the machinery is set up, it delivers Theorem \ref{01} with great ease, and
with a single treatment for all choices of the parameters $k$, $d$ and $n$.  
It transpires that the method should be widely applicable in related contexts. With a careful use of the ideas contained in this paper, the conclusions in Theorem \ref{01} for the equation \rf{05} can be refined to
$$  \mathrm{N}^*(B) = \frac{33-6\zeta(2)}{8\zeta(2)\zeta(3)} B\log B+ AB+O(B^{7/8}(\log B)^3)
$$
where $A$ is a suitable real number. We do not present the details here but refer to our forthcoming article
\cite{BB}. Another class of varieties where our methods can be applied is related to Vinogradov's mean value theorem. 
While the diagonal equation \rf{04} can be viewed as the ``germ'' that produces \rf{01} by inserting suitable products, we now consider a system
\be{08}  a_{l,0}y_0^l + a_{l,1}y_1^l + \ldots + a_{l,n} y_n^l =0 \quad (1\le l\le d) \ee
with integer coefficients $a_{l,j}$ as the germ. Then, the equations  
\be{09} \sum_{j=0}^n a_{l,j} (x_{1,j}x_{2,j} \cdots x_{k,j})^l = 0 \quad (1\le l\le d)
\ee
are homogeneous in the variables $\bfx_i = (x_{i,0},\ldots,x_{i,n})$,
and therefore the system defines another variety  in $\mathbb{P}^n(\QQ)^k$. As before, its rational points 
 correspond to solutions of \rf{09} in primitive vectors
$\bfx_i\in\ZZ^{n+1}$, with each $\bfx_i$ unique up to sign. The height of a rational point is now defined by
$ (|\bfx_1||\bfx_2|\cdots |\bfx_k|)^{n+1-\frac12 d(d+1)}$. Let $\mathrm{N}'(B)$ denote the number of rational points on the multiprojective variety defined by \rf{09} with height not exceeding $B$ and with $x_{i,j}\neq 0$ for $1\le i\le k$ and $0\le j\le n$. There is then a result for $\mathrm{N}'(B)$ that is very similar to Theorem \ref{Thm01}. Again, our method requires nearly optimal control on a certain moment of a suitable Weyl sum. 
Let $n_1(d)$ denote the smallest even natural number with the property that for any positive real number $\eps$ one has
$$ \int_{[0,1]^d} \Big| \sum_{1\le x\le P} {\mathrm e}^{2\pi{\mathrm i}(\al_1 x + \ldots+\al_d x^d)}\Big|^{n_1(d)}\,\dd\al \ll P^{n_1(d)-\frac12 d(d+1) +\eps}. $$

\begin{thm}\label{Thm02} Let $d,k,n$ be natural numbers with $n\ge n_1(d)$. Let $a_{l.j}$ $(1\le l\le d,\, 0\le j\le n)$ be non-zero integers. Then, there are a positive number $\delta$, a non-negative number $D$, and a monic polynomial $R\in\RR[X]$ of degree $k-1$ such that
$$\mathrm{N}' (B) = DBR(\log B) + O(B^{1-\delta}). $$
Moreover, the number $D$ is positive if and only if the system of equations \rf{08}
has non-trivial solutions in real numbers and in $p$-adic numbers, for all primes $p$.
\end{thm}

There is a widely held belief that $n_1(d) = d(d+1)$ should hold for all $d\in\NN$. This is trivial for $d=1$ and well-known for $d=2$ (see, for example, the discussion in \cite{BBMh}). Very recently Wooley \cite{W3} proved 
$n_1(3)=12$. With the arrival of efficient congruencing in epoque-making work of Wooley \cite{EC1}, early upper bounds on $n_1(d)$ due to Vinogradov and others remain of historical value. Progress with the new methods is still ongoing, but we now know that \cite{Wmulti}
$$ d(d+1) \le n_1(d) \le 2d(d-1)+2.$$ 

Our proof of Theorem \ref{Thm01} begins with counting primitive solutions to \rf{09} inside boxes 
$|\bfx_i|\le X_i$ $(1\le i\le k)$.  Irrespective of the shape of the box, connaisseurs of the circle method will find the task of supplying an asymptotic formula for this count  fairly routine. As a second step, one needs to express $\mathrm{N}(B)$ in terms of the count in boxes. More precisely, let $\bfm\in\NN^k$, and let $\theta(\bfm)$ denote the number of primitive $\bfx_i\in\ZZ^{n+1}$ satisfying \rf{09} and $ |\bfx_i|=m_i$, $x_{i,j}\neq 0$ $(1\le i\le k, \,0\le j\le n)$. Then, according to the opening paragraph of this memoir,
\be{010} \mathrm{N}(B) = 2^{-k} \sum_{m_1m_2\ldots m_k \le B^{1/(n+1-d)}} \theta(\bfm), \ee
and one desires to deduce an asymptotic formula for this sum from related formulae for 
\be{011}\Theta(X_1,\ldots,X_k)= \multsum{m_i\le X_i}{1\le i\le k} \theta(\bfm) \ee
that are the relevant sums for the count in boxes. The transition from \rf{011} to \rf{010} can be performed 
subject to very mild and general conditions, the actual definition of $\theta$ being irrelevant for this part of the argument. Hence, apart from the specific applications to problems such as those considered in Theorems \ref{Thm01} and \ref{Thm02}, we provide a method to evaluate divisor sums with weights similar to the one in 
\rf{010}. This tool should be of some utility in the analytic theory of numbers, and certainly constitutes the most novel feature of our analysis. Precise statements and a discussion of the underlying ideas of the 
transition process from \rf{011} to \rf{010} are deferred to next section that can be read independently of the rest of the paper. In Section 3 we study Weyl sums over products in preparation for the circle method work in Section 4. Each of these sections is equipped with a short introduction to the respective subject. In the short final section, Theorem \ref{Thm01} is deduced from the main conclusions formulated in Sections 2 and 4. Not only in structure but up to fine detail, the proof of Theorem \ref{Thm02} is very similar and therefore omitted.

\medskip

{\em Notation}. Owing to the wide range of topics discussed in this paper, 
a completely consistent notation would be in conflict with the various traditions involved. However, most of the notation used is standard. Small italics usually denote integers, letters from the middle of the alphabet like $k$, $n$ or $m$,
but also $d$
are natural numbers, and $j$ is a non-negative integer. The letter $p$ is reserved for primes. Capital italics $N,P,B,U,V,W,X,Y$ and $Z$ are used 
for the main parameters, and statements involving such a parameter are valid for its positive values exceeding $1$. Real variables are denoted by small greek letters, but also by $t$. These conventions apply with or without subscript whenever  the symbol does not obviously denote a function.

The letter $\eps$ always denotes a positive real number, and whenever $\eps$ occurs in a statement
it is asserted that the statement is valid for all positive numbers $\eps$. Constants implicit in Landau or Vinogradov symbols may depend on $\eps$. Note that this convention allows us to conclude from $A\ll P^\eps$ and $B\ll P^\eps$ that $AB\ll P^\eps$, for example.

We make use of vector notation frequently, with some uncommon elements now to be described. 
Boldface denotes vectors, with entries written in the corresponding italic font, and the dimension
may vary from one occurrence to the next. When $\bfx=(x_1,\ldots,x_r)$, $\bfy=(y_1,\ldots,y_s)$ and
$\bfz=(z_1,\ldots,z_t)$, we  use $(\bfx,\bfy,\bfz)$ 
to denote the $r+s+t$-tuple $(x_1,\ldots,x_r,y_1,\ldots,y_s,z_1,\ldots,y_t)$.
Also, we will have to permute the entries of a vector. Let $S_k$
denote the symmetric group on $k$ elements. Then, for $\sgm\in S_k$ and
$\bfx\in\RR^k$, write
$$ \mbox{}^\sgm\bfx = (x_{\sgm(1)},\ldots,x_{\sgm(k)}). $$
For  functions $H$ defined on $\NN^k$ or $[1,\infty)^k$, we define
$$ H_\sgm (\bfx) = H (\mbox{}^\sgm\bfx). $$ 
Inequalities between vectors are to be interpreted as the system of 
inequalities given by the components. Thus,  for $\bfx, \bfX\in \RR^k$, the 
system of $k$ inequalities
$x_j\le X_j$ $(1\le j\le k)$ is abbreviated to $\bfx\le \bfX$. Whenever
$X_j\ge 1$ for all $j$ with $1\le j\le k$, we write
$$  \langle \bfX\rangle = X_1X_2\ldots X_k. $$   

The number of divisors of $n$ is denoted by $\tau(n)$, Euler's totient function is $\varphi(n)$,
and the M\"obius function is $\mu(n)$. The highest common 
factor of $a$ and $b$ is $(a;b)$. We put $e(\al)=\exp(2\pi\mathrm{i}\al)$, and denote Riemann's zeta function by $\zeta(s)$.

\section{The hyperbola method}

\subsection{The transition theorem}
Not only a few problems in the theory of numbers depend implicitly or explicitly
on the asymptotic evaluation of the sum
\be{211} \Ups(N) = \sum_{u_1u_2\ldots u_k\le N} h(\bfu). \ee
Here the dimension $k$ and the arithmetical function $h: \NN^k \to \CC$ are
intrinsic to the application 
at hand. Perhaps the most familiar cases are the divisor 
problems of Dirichlet and Piltz where one chooses $h(\bfu)=1$. As we have
indicated in the introductory section of the present communication, the
counting problems discussed in Theorems \ref{Thm01} and \ref{Thm02} also reduce to 
sums of the type \rf{211}. 
Yet, a successful treatment of
the cognate box sums
\be{212} H(X_1,\ldots,X_k) = \multsum{1\le x_j\le X_j}{1\le j\le k} h(\bfx)\ee
is often easier, in particular  in cases
where $h(\bfu)$ is the number 
of solutions of a certain diophantine system. Typically, the condition on the 
product of the variables $u_1u_2\ldots u_k$ in \rf{211} will be in conflict 
with a direct
use of the Hardy-Littlewood method for the diophantine problem at hand.

Whenever the evaluation of the box sums is within the compass of existing technology,
one is led to the question whether suitable 
asymptotic formulae for \rf{212} contain sufficient 
information to deduce an allied
formula for $\Ups(N)$. We shall provide an affirmative answer when $h$ 
takes real non-negative values only, and when the leading term in the 
asymptotic formula for $H(\bfX)$ is a pure power of $X_1X_2\ldots X_k$. 
These conditions are not infrequently met in practice, and Theorem \ref{S2.1}
below should be useful in areas other than those discussed in this
paper. However, it should be noted that the conditions on $h$ 
formulated in the preamble to Theorem \ref{S2.1} have been tuned for our 
immediate needs. The underlying arguments work in  broader generality. 
For example, if the leading term in an asmptotic formula for $H(\bfX)$ contains
logarithms, these may be accommodated by the method now to be described. For another development of our ideas,
see Schindler \cite{Sch}.

We begin by introducing the class of functions $h$ to which our theory applies.
Fix positive real numbers $\al, c,\del$ with $\del<\min(1,\al)$.   
A set $\cal H$ of 
arithmetical functions $h: \NN^k\to [0,\infty)$ will be referred to as a
{\em family satisfying condition (I) with respect to $(\al, c,\del)$} if the following 
holds:

\medskip\noindent
(I) {\em For any $h\in\cal H$ there is a real number $c_h\in [0,c]$
such that  the asymptotic formula
$$ \sum_{\bfx \le \bfX} h(\bfx) = c_h \langle\bfX\rangle^\al+ 
O(\langle\bfX\rangle^\al(\min_{1\le j\le k} X_j)^{-\del}) $$
holds uniformly in $h\in\cal H$ and $X_j\ge 1$ $(1\le j\le k)$.}

\medskip\noindent
Now consider a family  $\cal H$ satisfying 
condition (I) with respect to $(\al, c,\del)$. Further,  suppose that the real numbers $\nu$  
and $D$ satisfy 
$ 0<\nu\le 1$ and $ D\ge 0 $. 
The set  $\cal H$ is called an
$(\al,c,D,\nu,\del)$-{\em family}, provided the two further conditions (II) and (III)
are satisfied:

\medskip\noindent
(II) {\em For $h\in\cal H$ and $r\in\NN$ with $1\le r\le k-1$, there
exists an arithmetical function $c_{h,r}: \NN^r \to [0,\infty)$
such that  for any $\bfu\in\NN^r$ the asymptotic formulae 
$$ \sum_{\bfv\le \bfV}  h(\bfu,\bfv) = c_{h,r}(\bfu) 
\langle\bfV\rangle^\al 
+ O(\langle\bfV\rangle^\al |\bfu|^D
(\min V_j)^{-\del}) $$
hold uniformly for $h\in\cal H$, $V_j \ge 1 $ and 
$|\bfu|\le\langle\bfV\rangle^\nu$, }

\medskip\noindent
(III) {\em For all $h\in\cal H$ and $\sgm\in S_k$, one has $h_\sgm\in\cal H$.}
\bigskip

Note that the condition (I) is symmetric with respect to the 
indices $1\le j\le k$ whereas (II) alone is not. However, by (III),
one may apply (II) to $h_\sgm$, for any $\sgm\in S_k$.
Hence, one may choose any $r$ indices and sum $h(\bfx)$ over the 
corresponding variables. There is then an asymptotic formula for this
sum, similar to the one in (II). Notice also that (I) can be interpreted
as the case $r=0$ of (II).
\medskip

We are ready to announce the principal result of this chapter.

\begin{thm}\label{S2.1}  Let $k\ge 2$, and let $\cal H$ 
be an $(\al, c,D,\nu,\del)$-family of arithmetical
functions $h:\NN^k \to [0,\infty)$. For any $h\in\cal H$, let $\Ups(N)$
be defined by \rf{211}. There exists a positive number $\eta$
with the property that for any $h\in\cal H$ there is  a 
polynomial $P_h\in\RR [x]$ of degree at most $k-2$
such that the asymptotic formula
\be{T} \Ups(N) = 
 N^\al \Big(\frac{c_h \al^{k-1}}{(k-1)!} (\log N)^{k-1}  +      P_h(\log N)\Big) + O(N^{\al-\eta})\ee
holds uniformly in $h\in\cal H$. 
\end{thm}

The proof of Theorem \ref{S2.1} produces an explicit value for $\eta$
in terms of $k,\al, c,D,\nu$ and $\del$, but it will not be very large, as
part of the argument is based on induction on $k$. Therefore no attempt has been made
to record the optimal value for $\eta$ that our methods could establish.

In the next section, we begin by 
observing a certain rigidity within the functions $c_{h,r}$ in (II). Then,
in the following three sections, we collect several estimates of 
preparatory character. Some of the asymptotic relations obtained here
may be of independent interest. We highlight the
light-weight version of Theorem \ref{S2.1} in Theorem \ref{S2.7} below. The main
argument leading to a proof of Theorem \ref{S2.1} is presented
in the closing section. 

\subsection{Families of arithmetical functions}
We begin our discussion with a rough yet useful estimate.

\begin{lem}\label{Lextra} Let  $\cal H$ denote a family satisfying 
condition (I) with respect to $(\al, c,\del)$. Then, uniformly in $h\in \cal H$, one has 
$h(\bfu) \ll \langle \bfu \rangle^\al$. 
\end{lem}

\noindent {\em Proof}.
This is immediate  from (I), on taking $\bfX=\bfu$. 

\medskip

The principal observation in this section is that 
the functions $c_{h,l}$ associated with an  $(\al, c,D,\nu,\del)$-{\em family} 
form another such family. The following lemma makes this precise.

\begin{lem}\label{L2.2}
 { Let $\cal H$ be an  $(\al, c,D,\nu,\del)$-{\em family}
of arithmetical
functions $h:\NN^k \to [0,\infty)$. Then, for any $1\le l\le k-1$, the set of
functions $c_{h,l}: \NN^l \to [0,\infty)$ with $h\in\cal H$ forms an
 $(\al, c,D,1,\del)$-{family}, and one has}
\be{221}
\sum_{\bfy\le\bfY}  c_{h,l}(\bfy) = c_h \langle\bfY\rangle^\al 
+ O(\langle\bfY\rangle^\al(\min Y_j)^{-\del}).
\ee 
\end{lem}

\noindent {\em Proof}.
We lauch the proof with the demonstration of \rf{221}. In the interest of
brevity, write $m=k-l$. Let $Z\ge 1$, and put $\bfZ= (Z,\ldots,Z)\in\RR^m$.
Then, as a special case of (II), the asymptotic relation
\be{222} \sum_{\bfz\le \bfZ} h(\bfy,\bfz) = c_{h,l}(\bfy)Z^{\al m} 
+ O(Z^{\al m-\del}|\bfy|^D)  \ee
holds uniformly for all $\bfy\in\NN^l$ with $|\bfy|\le Z^{m\nu}$. This may be 
summed over a box. Note that $m\ge 1$, so that whenever 
$\bfY\in\RR^l$ satisfies $|\bfY|\le Z^\nu$, one certainly has
$$ \sum_{\bfy\le\bfY}  \sum_{\bfz\le \bfZ} h(\bfy,\bfz)
= Z^{\al m} \sum_{\bfy\le\bfY}  c_{h,l}(\bfy) + O(\langle\bfY\rangle^{D+1} Z^{\al m-\del}). $$
On the other hand, the sum on the left may be evaluated by (I).
For $Z\ge |\bfY|$, this yields
$$ \sum_{\bfy\le\bfY}  \sum_{\bfz\le \bfZ} h(\bfy,\bfz) 
= c_h \langle\bfY\rangle^\al Z^{\al m} +  O( Z^{\al m} \langle\bfY\rangle^\al(\min Y_j)^{-\del}).$$
Divide by $Z^{\al m}$ and recall that $\nu\le 1$ to conclude from the last two formulae that for 
$Z\ge|\bfY|^{1/\nu}$ one has
$$ \sum_{\bfy\le\bfY}  c_{h,l}(\bfy) = c_h \langle\bfY\rangle^\al 
+ O(\langle\bfY\rangle^\al(\min Y_j)^{-\del}) + O(\langle\bfY\rangle^{D+1} Z^{-\del}). $$
With $Z\to\infty$, the asymptotic formula \rf{221} follows. In particular,
this verifies condition (I) for the family $c_{h,l}$. Note that the implicit
constant in \rf{221} is inherited from the 
conditions (I) and (II) for the family $\cal H$, and is therefore 
uniform in $h$.

Next, we establish (II). It will be appropriate to adopt the notation
from the previous argument. In addition, let $1\le r\le l-1$ and
$\bfu\in\NN^r$, $\bfv\in\NN^{l-r}$. We choose $\bfy=(\bfu,\bfv)$ in \rf{222}
and sum over a box for $\bfv$. Then, provided that $|\bfu|\le Z^{m\nu}$,
$|\bfV|\le Z^{m\nu}$, one finds that
$$ \sum_{\bfv\le\bfV}  \sum_{\bfz\le \bfZ} h(\bfu,\bfv,\bfz)
= Z^{\al m} \sum_{\bfv\le \bfV} c_{h,l}(\bfu,\bfv) + O(|\bfu|^D\langle\bfV\rangle^{D+1} Z^{\al m-\del}).
$$ 
Alternatively, one may use (II) to evaluate the left hand side above.
When $|\bfV|\le Z$ and $|\bfu| \le \langle\bfV\rangle^\nu Z^{m\nu}$, this yields
$$ \sum_{\bfv\le\bfV}  \sum_{\bfz\le \bfZ} h(\bfu,\bfv,\bfz)
= c_{h,r}(\bfu) \langle\bfV\rangle^\al Z^{\al m} + O(|\bfu|^D\langle\bfV\rangle^\al Z^{\al m} (\min V_j)^{-\del}).$$
One may now proceed as before: a comparison of the last two displays 
delivers the preliminary estimate
$$ \sum_{\bfv\le \bfV} c_{h,l}(\bfu,\bfv) =  c_{h,r}(\bfu) \langle\bfV\rangle^\al
 + O(|\bfu|^D\langle\bfV\rangle^\al  (\min V_j)^{-\del}) +  
O(|\bfu|^D\langle\bfV\rangle^{D+1} Z^{-\del}), $$
subject to the lower bounds  on $Z$ mentioned earlier. With $Z\to\infty$,
the term on the far right disappears, and any $\bfu$ satisfies $|\bfu| \le  Z^{m\nu}$
when $Z$ is large enough. This gives
$$ \sum_{\bfv\le \bfV} c_{h,l}(\bfu,\bfv) =  c_{h,r}(\bfu) \langle\bfV\rangle^\al
 + O(|\bfu|^D\langle\bfV\rangle^\al  (\min V_j)^{-\del}) $$
uniformly in $\bfu$. This proves (II) for $c_{h,l}$ in place of $h$, and with $\nu=1$. Moreover, the last 
estimate is uniform with respect to $h\in\cal H$ because the 
implicit constant can be traced back to the one in (II). 
Note the recurrent appearance
of $c_{h,r}$ on the right hand side.

Finally, we have to check the condition (III). Let $\pi\in S_l$
be a permutation of $\{1,\ldots,l\}$, and let $\sgm\in S_k$ be a
permutation of $\{1,\ldots,k\}$ with $\sgm(j)=\pi(j)$ for $1\le j\le l$.
Then, by (II), it is immediate that $(c_{h,l})_\pi = c_{h_\sgm,l}$, which
confirms (III), completing the proof of the lemma.

\bigskip

Let $\cal H$ be an $(\al,c,D,\nu,\del)$-family of arithmetical functions 
defined on $\NN^k$, and let $1\le l\le k-1$. For $h\in\cal H$, one may fix
$l$ of the variables and consider $h$ as a function on $\NN^{k-l}$. This 
process yields another family of functions with similar properties. More
precisely, choose a real number $A\ge \al$ and consider, for any $h\in\cal H$
and $\bfw\in\NN^l$, the function
$$ g=g_{h,\bfw}: \NN^{k-l} \to [0,\infty), \quad \bfy \mapsto \langle\bfw
\rangle^{-A} h(\bfw,\bfy). $$
Let ${\cal H}_l$ denote the set of all these functions. Since $\cal H$ 
satisfies (III), the same is true of ${\cal H}_l$. 

\begin{lem}\label{L2.3}
{ Let $\cal H$ and $l$ be as above, and suppose that
$A\ge D+ (k+1)\al + \nu^{-1}(1+\al)$. 
Then, for sufficiently large $c'$, the set  ${\cal H}_l$ is 
an $(\al,c',D, \nu,\del)$-family.}
\end{lem}

\noindent
{\em Proof}. By Lemma \ref{L2.2}, we may apply Lemma \ref{Lextra}
to $c_{h,l}$ in place of $h$.  This yields the bound 
\be{chl}\langle\bfw\rangle^{-\al} c_{h,l}(\bfw)\ll 1. \ee
We
 now proceed to show that uniformly for $h\in\cal H$ and
$\bfw\in\NN^l$, one has
\be{I1} 
 \sum_{\bfy\le\bfY} \f{h(\bfw,\bfy)}{\langle \bfw \rangle^A}
= \f{c_{h,l}(\bfw)}{ \langle \bfw \rangle^A} \langle \bfY \rangle^\al
+ O (\langle \bfY \rangle^\al (\min Y_j)^{-\del}).
\ee
Once this is established, one obtains (I) for $g_{h,\bfw}$, with
$\langle\bfw\rangle^{-A} c_{h,l}(\bfw)$ in the role of $c_h$, and since
$A\ge \al$, 
the existence of $c'$ is a consequence of \rf{chl}.

For $|\bfw|\le \langle\bfY\rangle^\nu$, one notes that $A\ge D$ to realize that
 the asymptotic formula \rf{I1}
is a weakened form of (II). Thus, it suffices to confirm \rf{I1} 
when  $|\bfw| > \langle\bfY\rangle^\nu$. In this case, we apply \rf{chl}
to see that
$$ \f{c_{h,l}(\bfw)}{\langle\bfw\rangle^A} \langle\bfY\rangle^\al
\ll  \langle\bfw\rangle^{\al-A} \langle\bfY\rangle^\al \ll 1.
$$
Similarly, by Lemma \ref{Lextra}, 
$$ \sum_{\bfy\le \bfY} \f{h(\bfw,\bfy)}{\langle\bfw\rangle^A}
\ll \langle\bfw\rangle^{\al-A} \langle\bfY\rangle^{\al+1} \ll 1.$$
Since $\del\le \min(1,\al)$, one concludes that the two explicit terms in
\rf{I1} are both bounded above by $O(\langle\bfY\rangle^\al (\min Y_j)^{-\del})$.
In particular, this proves \rf{I1}.

The proof of (II) for ${\cal H}_l$ is very similar. Let $1\le r\le k-l-1$, and
$\bfv\in\NN^r$. We need to confirm that whenever $|\bfv
|\le \langle\bfZ\rangle^\nu$, one has
\be{II1} \sum_{\bfz\le\bfZ} \f{h(\bfw,\bfv,\bfz)}{\langle\bfw\rangle^A}
= \f{ c_{h,l+r}(\bfw,\bfv)}{\langle\bfw\rangle^A} \langle\bfZ\rangle^\al
+ O(\langle\bfZ\rangle^\al (\min Z_j)^{-\del} |\bfv|^D).
\ee
As before, we begin with the case where $|\bfw|\le \langle\bfZ\rangle^\nu$.
Here, it suffices to recall that $A\ge D$ which implies that
$|(\bfw,\bfv)|^D\langle\bfw\rangle^{-A}\le |\bfv|^D$. Hence, in this case
\rf{II1}
 is immediate from
(II). 

When $|\bfw| > \langle\bfZ\rangle^\nu$, one uses Lemma \ref{Lextra} to 
deduce that
$$ \sum_{\bfz\le \bfZ} \f{h(\bfw,\bfv,\bfz)}{\langle\bfw\rangle^A}
\ll \langle\bfw\rangle^{\al-A} \langle\bfv\rangle^{\al}
\langle\bfZ\rangle^{\al+1} \ll 1,$$
and a similar but simpler estimation
gives
$$ \f{ c_{h,l+r}(\bfw,\bfv)}{\langle\bfw\rangle^A} \langle\bfZ\rangle^\al
\ll 1. $$
Now \rf{II1} follows as above. 

We have already remarked that ${\cal H}_l$ satisfies (III).
The proof of Lemma \ref{L2.3} is complete.

\subsection{Preparatory lemmata} Before we may announce our first estimate,
we must introduce some more notation. In order to avoid repetitious comments
concerning uniformity of implicit constants, we remark that for the remainder of this 
chapter,
these constants will depend on an auxiliary non-negative integer $j$, the dimension $k$, the family $\cal H$,
and on $\eps$ where appropriate. In particular, the constants depend on
the parameters $\alpha, c,D, \nu, \delta$.
It is relevant to note,
however, that these constants will not depend on the individual $h\in\cal H$.

Now let  
 $$ \Delta^{(k)} = \{\mathbf{t} \in \Bbb{R}^k : 
1 < t_1 <  t_2 < \ldots < t_k\}, \quad
\Delta^{(k)}(X) =  \{\mathbf{t} \in \Delta^{(k)}: t_k\le X\}.$$
The next lemma evaluates the integral
$$ I_{h,j}(X) = \int_{\Delta^{(k)}(X)} 
\frac{(\log\langle\mathbf{t}\rangle)^j}{\langle\mathbf{t}\rangle^{\alpha+1}} 
\sum_{\mathbf{x} \leq \mathbf{t}} h(\mathbf{x}) \,\dd\mathbf{t}. $$

\begin{lem}\label{L2.4} {Let $\mathcal{H}$ be 
an $(\alpha, c,D, \nu, \delta)$-family and $j\in\NN_0$. Then there exists 
a positive number 
$\eta$ such that for all 
$h \in \mathcal{H}$ and suitable real polynomials $q_{h,j}$ one has
$$ I_{h,j}(X) = q_{h,j}(\log X) + O(X^{-\eta}).$$
}
\end{lem}

It would be possible to compute the degree and the leading coefficient of $q_{h,j}$, 
but this will not be relevant later.

\medskip
\noindent
{\em Proof.}  For $k=1$, we may rewrite (I) in the form
\begin{displaymath}
  \sum_{x \leq t} h(x) = c_ht^{\alpha} + E(t)
\end{displaymath}
where $E(t)=E_h(t)$ is piecewise continuous and satisfies 
$E(t) \ll t^{\alpha-\delta}$. Hence, by straightforward estimates, one finds that
$$
  \int_1^X \frac{(\log t)^j}{t^{\alpha+1}} \sum_{x \leq t} h(x)\,\dd t = \frac{c_h}{j+1} (\log X)^{j+1} + D_{h,j} + O(X^{-\delta} (\log X)^j) $$
in which the real number $ D_{h,j} $ is given by the convergent integral
$$ D_{h,j}= \int_1^{\infty} E(t)t^{-\alpha-1} (\log t)^j \,\dd t.
$$
This confirms the assertion of the lemma when $k=1$. More precisely, in this case,
any $0<\eta<\del$ is admissible.

We proceed by induction and suppose that
$k>1$, and that the lemma has been established for all smaller values of $k$. 
We split the set $ \Delta^{(k)}$ into $k$ disjoint subsets.  
To describe this dissection, put $t_0=1$ and
\be{defb}  \beta = \min (\nu, \del (2D+4k)^{-1}).
\ee
For $0\le l < k$, let 
$$
  \Delta^{(k, l)} = \{\mathbf{t} \in  \Delta^{(k)} :  
t_i> t_{i+1}^\beta\;(l<i<k),\; t_l\le t_{l+1}^\beta\} $$
and write
$ \Delta^{(k, l)}(X)= 
\Delta^{(k, l)}\cap \Delta^{(k)}(X).
$
Then indeed $\Delta^{(k)}(X)$ is the disjoint union of $ \Delta^{(k, l)}(X)$ 
with $0\le l<k$,
so that we now have
\be{232} I_{h,j}(X) = \sum_{l=0}^{k-1} I_{h,j,l}(X)  \ee
where 
\be{230}I_{h,j,l}(X) = \int_{\Delta^{(k,l)}(X)} 
\frac{(\log\langle\mathbf{t}\rangle)^j}{\langle\mathbf{t}\rangle^{\alpha+1}} 
\sum_{\mathbf{x} \leq \mathbf{t}} h(\mathbf{x}) \,\dd\mathbf{t}. 
\ee

The integral $I_{h,j,0}(X)$ can be computed in much the same way as in the 
treatment of the case $k=1$. For  $\mathbf{t} \in \Delta^{(k, 0)}(X)$, we define
the function $E_h(\mathbf{t})$ via
\be{233}
\sum_{\mathbf{x} \leq \mathbf{t}} h(\mathbf{x}) = c_h\langle\mathbf{t}\rangle^{\alpha} 
+ E_h(\mathbf{t}).
\ee
Then, by (I), one has  $E_h(\mathbf{t}) \ll \langle\mathbf{t}\rangle^{\alpha} t_1^{-\delta}$.
The simple bound $ \log\langle\mathbf{t}\rangle \le k \log |\mathbf{t}|$ 
suffices to confirm
that for $Z\ge 1$ one has
\begin{eqnarray*}
\lefteqn{
\int_{\Delta^{(k, 0)}(2Z)\setminus \Delta^{(k, 0)}(Z)} 
\frac{(\log\langle\mathbf{t}\rangle)^j}{\langle\mathbf{t}\rangle^{\alpha+1}} |E_h(\mathbf{t})|
\,\dd\mathbf{t} }
\\
& \ll & \int_Z^{2Z} \int_{t_k^{\beta}}^{t_k}  \int_{t_{k-1}^{\beta}}^{t_{k-1}} \cdots 
\int_{t_2^{\beta}}^{t_2} \langle\mathbf{t}\rangle^{-1} t_1^{-\delta} ( \log t_k)^j
\,\dd t_1 \ldots \dd t_k
\ll Z^{-\del\beta^k}. 
\end{eqnarray*}
On summing over dyadic ranges, it follows that the integral
$$  D_{h,j,0}  = \int_{\Delta^{(k, 0)}} 
\frac{(\log\langle\mathbf{t}\rangle)^j}{\langle\mathbf{t}\rangle^{\alpha+1}} E_h(\mathbf{t})
\,\dd\mathbf{t} $$
exists, and that this differs from the same integral over $\Delta^{(k, 0)}(X)$
by $O(X^{-\del\beta^k})$. Consequently, by \rf{230} and \rf{233},
\be{234} I_{h,j,0}(X) =  c_h\int_{\Delta^{(k, 0)}(X)} 
 \frac{(\log\langle\mathbf{t}\rangle)^j}{\langle\mathbf{t}\rangle} \,\dd\mathbf{t}  + 
D_{h,j,0} 
+ O(X^{- \delta\beta^k}).
\ee

At this point, we interrupt the treatment of $ I_{h,j,0}(X)$ and turn to
$  I_{h,j,l}(X)$ with $1\le l<k$. 
 For $\mathbf{t}\in\RR^k$, 
write $\mathbf{t} = (\mathbf{t}', \mathbf{t}'')$ 
with $\mathbf{t}' = (t_1, \ldots, t_l)$ and 
$\mathbf{t}'' = (t_{l+1}, \ldots, t_k)$. 
By \rf{defb}, one has  $|\mathbf{t}'| \le |\mathbf{t}''|^\nu$
 for all $\mathbf{t}
\in \Delta^{(k, l)}$. Hence,
by (II), the function $E_{h,l}(\mathbf{x}', \mathbf{t}'')$ defined through
the equation
\be{235}
  \sum_{\mathbf{x}'' \leq \mathbf{t}''} h(\mathbf{x}', \mathbf{x}'') = 
c_{h,l}\langle\mathbf{x}') \langle\mathbf{t}''\rangle^{\alpha} + E_{h,l}(\mathbf{x}', \mathbf{t}'')
\ee
satisfies the bound $E_{h,l}(\mathbf{x}', \mathbf{t}'')
\ll  \langle\mathbf{t}''\rangle^{\alpha} t_{l+1}^{-\delta}|\mathbf{x}'|^D$ 
uniformly in $\mathbf{x}' \leq \mathbf{t}'$. The strategy is now the same as
in the treatment of $I_{h,j,0}(X)$. We shall sum \rf{235} over
 $\mathbf{x}'\le\mathbf{t}' $ and insert the result into \rf{230}. With this 
end in view, put
$$ R_l (\mathbf{t}) = \sum_{\mathbf{x}'\le\mathbf{t}'} 
E_{h,l}(\mathbf{x}', \mathbf{t}''). $$
This defines a function on $ \Delta^{(k, l)}$ satisfying
$ R_l (\mathbf{t})\ll 
\langle\mathbf{t}'\rangle  \langle\mathbf{t}''\rangle^{\alpha} t_{l+1}^{-\delta}t_l^D$.
For $Z\ge 1$ we now see that 
\begin{eqnarray}\nonumber
\lefteqn{
\int_{\Delta^{(k, l)}(2Z)\setminus \Delta^{(k, l)}(Z)} 
\frac{(\log\langle\mathbf{t}\rangle)^j}{\langle\mathbf{t}\rangle^{\alpha+1}} 
|R_l(\mathbf{t})|\,\dd\mathbf{t} }
\\
& \ll & \int_{{\cal R}_l(Z)} \langle\mathbf{t}''\rangle^{-1} t_{l+1}^{-\delta} (\log t_k)^j
\int_{\Delta^{(l)}(t_{l+1}^{\beta})} \langle\mathbf{t}'\rangle^{-\al}t_l^D \,\dd\mathbf{t}'\,\dd\mathbf{t}''
\label{236}\end{eqnarray}
where
$$ {\cal R}_l(Z) = \{{\mathbf t}'': Z<t_k\le 2Z,\; t_{i+1}^{\beta} \le t_i\le t_{i+1} \;(l<i<k)\}. $$

We bound the inner integral by brute force. First one notes that
$ \langle\mathbf{t}'\rangle^{-\al} \le 1$, and consequently that 
$ \langle\mathbf{t}'\rangle^{-\al}t_l^D \le t_{l+1}^{D\beta}$ holds for all 
${\mathbf t}'\in\Delta^{(l)}(t_{l+1}^{\beta})$. Then, since the measure of 
$\Delta^{(l)}(t_{l+1}^{\beta})$ does not exceed $t_{l+1}^{l\beta}$, we find that the 
inner integral in \rf{236} is bounded above by $t_{l+1}^{(D+k)\beta}$, irrespective of the actual value of $l$. By \rf{defb}, we have $(D+k)\beta < \f12\del$. Hence, by \rf{236},
$$
\int_{\Delta^{(k, l)}(2Z)\setminus \Delta^{(k, l)}(Z)} 
\frac{(\log\langle\mathbf{t}\rangle)^j}{\langle\mathbf{t}\rangle^{\alpha+1}} 
|R_l(\mathbf{t})|\,\dd\mathbf{t} 
 \ll  \int_{{\cal R}_l(Z)} \langle\mathbf{t}''\rangle^{-1} t_{l+1}^{-\delta/2}\,
\dd\mathbf{t}'' \ll Z^{-\del\beta^k/2}.  $$
Much as before, on summing over dyadic ranges, it follows that
the integral
$$ D_{h,j,l} = \int_{\Delta^{(k, l)}} 
\frac{(\log\langle\mathbf{t}\rangle)^j}{\langle\mathbf{t}\rangle^{\alpha+1}} R_l(\mathbf{t})
\,\dd\mathbf{t} $$
exists and differs from the same integral over ${\Delta^{(k, l)}}(X)$
by $O(X^{-\del\beta^k/2)})$. By \rf{235} and \rf{230}, we infer that
\be{236a}
I_{h,j,l}(X) = \int_{\Delta^{(k, l)}(X)} 
\frac{(\log\langle\mathbf{t}\rangle)^j}{\langle\mathbf{t}''\rangle\langle\mathbf{t}'\rangle^{\alpha+1}}
 \sum_{\mathbf{x}' \leq \mathbf{t}'}c_{h,l}(\mathbf{x}')\,\dd\mathbf{t}
+D_{h,j,l}+O(X^{-\del\beta^k/2}). 
\ee
Let $K_{h,j,l}(X)$ denote the integral on the right hand side of the previous display.
Since $\log \langle\mathbf{t}\rangle = \log \langle\mathbf{t}'\rangle + \log \langle\mathbf{t}''\rangle$,
binomial expansion shows that
\be{236b}
K_{h,j,l}(X) = \sum_{r =0}^j \left(\begin{matrix} j\\r\end{matrix} \right)
\int_{{\cal S}_l(X)}  \frac{(\log\langle\mathbf{t}''\rangle)^{j-r} }{\langle\mathbf{t}''\rangle}
\int_{\Delta^{(l)}(t_{l+1}^{\beta})} 
   \frac{(\log\langle\mathbf{t}'\rangle)^{r} }{\langle\mathbf{t}'\rangle^{\alpha+1}}\sum_{\mathbf{x}' \leq \mathbf{t}'} c_{h,l}(\mathbf{x}')\, \dd\mathbf{t}' \,\dd\mathbf{t}'' 
\ee
in which
$$ {\cal S}_l(X) = \{\mathbf{t}'': t_k\le X,\; t_{i+1}^{\beta}\le t_i \le t_{i+1}\;
(l<i<k)\}. $$

To compute the inner integral in \rf{236b}, we apply the induction hypothesis together
with Lemma \ref{L2.2}. It then follows that this integral equals 
$Q_{h,l,r}(\log t_{l+1}^\beta) + F_{h,l,r}(t_{l+1}^{\beta})$, with  a suitable real
polynomial $Q_{h,l,r}$
and a suitable function $ F_{h,l,r}$ satisfying the estimate
$ F_{h,l,r}(t) \ll t^{-\eta(l)}$. Here, $\eta(l)$ is the positive number that the induction
hypothesis produces for the family $c_{h,l}$. 
An argument similar to the one used 
around \rf{236} and leading to  \rf{236a} now shows that for there is a real number 
$E_{h,j,l,r}$ such that
$$ \int_{{\cal S}_l(X)}  \frac{(\log\langle\mathbf{t}'' \rangle)^{j-r} }{\langle\mathbf{t}''\rangle}
F_{h,l,r}(t_{l+1}^{\beta}) \,\dd\mathbf{t}'' = E_{h,j,l,r} + O(X^{-\eta(l)\beta^k/2}).$$
From \rf{236b} we  now deduce that there is a real number, say $E_{h,j,l}$,
such that
$$ K_{h,j,l}(X)= E_{h,j,l} + \int_{{\cal S}_l(X)} \sum_{r =0}^j \left(\begin{matrix} j\\r\end{matrix} \right)
 \frac{(\log\langle\mathbf{t}''\rangle)^{j-r} }{\langle\mathbf{t}''\rangle}
Q_{h,l,r}(\log t_{l+1}^{\beta}) \,\dd\mathbf{t}''+ O(X^{-\eta(l)\beta^k/2}). $$
Hence, whenever $0<\eta <\f12 \beta^k \min (\eta(l), \delta)$, the formula \rf{236a}
yields
$$ I_{h,j,l}(X)= E_{h,j,l} + \int_{{\cal S}_l(X)} \sum_{r =0}^j \left(\begin{matrix} j\\r\end{matrix} \right)
 \frac{(\log\langle\mathbf{t}''\rangle)^{j-r} }{\langle\mathbf{t}''\rangle}
Q_{h,l,r}(\log t_{l+1}^{\beta}) \,\dd\mathbf{t}''+ O(X^{-\eta}). $$ 

To compute the remaining integral here, 
one expands $Q_{h,l,r}(\log t_{l+1}^{\beta})(\log\langle\mathbf{t}''\rangle)^{j-r}$
as a polynomial in $\log t_i$, $l<i\le k$. We may then rewrite this
integral as a linear combination of integrals of the type
\begin{eqnarray*}\lefteqn{  \int_{{\cal S}_l(X)} \frac{(\log t_k)^{b_k} \cdots (\log t_{l+1})^{b_{l+1}}}{t_k \cdots t_{l+1}} \,\dd\mathbf{t}'' }
\\ &=&
  \int_1^X \int_{t_k^{\beta}}^{t_k} \cdots \int_{t_{l+2}^{\beta}}^{t_{l+2}}  
\frac{(\log t_k)^{b_k} \cdots (\log t_{l+1})^{b_{l+1}}}{t_k \cdots t_{l+1}}   
\,\dd t_{l+1} \cdots \dd t_k
\end{eqnarray*}
with $b_{i+1}, \ldots, b_k \in \Bbb{N}_0$. These integrals can be computed explicitly,
and are polynomials in $\log X$. Similarly, the integral in \rf{234}
is a polynomial in $\log X$. This shows that all $I_{h,j,l}(X)$ are polynomials in
$\log X$, up to an error not exceeding $O(X^{-\eta})$. By \rf{232},
the same is then true for $I_{h,j}(X)$, completing the induction.

\begin{lem}
\label{L2.6} {Let $\mathcal{H}$ be 
an $(\alpha, c,D, \nu, \delta)$-family and  $j\in\NN_0$. Then there exists a positive number
$\eta$ such that for all $h\in\cal H$ and suitable real polynomials $Q_{h,j}$ one has
$$  \int_{[1,X]^k}\frac{(\log\langle\mathbf{t}\rangle)^j}{\langle\mathbf{t}\rangle^{\alpha+1}} 
\sum_{\mathbf{x} \leq \mathbf{t}} h(\mathbf{x}) \,\dd\mathbf{t}   = Q_{h,j}(\log X) + O(X^{-\eta}).$$
}
\end{lem}

For a proof, we only need to observe that 
for any $\mathbf{t}\in [1,X]^k$ satisfying $t_i\neq t_l$ for all $i\neq l$
there is a unique $\sgm\in S_k$ such that 
$\mbox{}^\sgm\mathbf{t}\in\Delta^{(k)}(X)$. Since the set of all 
$\mathbf{t} \in \Bbb{R}^k$ where at least two coordinates are equal is a set
of measure $0$, we may conclude that
$$
 \int_{[1, X]^k} \frac{(\log\langle\mathbf{t}\rangle)^j}{\langle\mathbf{t}\rangle^{\alpha+1}} \sum_{\mathbf{x} \leq \mathbf{t}} h(\mathbf{x}) \,\dd\mathbf{t}  = 
\sum_{\sgm\in S_k} I_{h_\sgm,j}(X),
$$
and the lemma follows from Lemma \ref{L2.4}.

\subsection{A mean value estimate}
Our  next result  is a light-weight version of Theorem \ref{S2.1}. 
It features the real number
\begin{displaymath}
  V_{k, j} = \int_{[0, 1]^k} (\xi_1 + \ldots + \xi_k)^j \,\dd\bm \xi, 
\end{displaymath}
defined whenever $k\in\NN$, $j\in\NN_0$. One may calculate this integral
elementarily to obtain the alternative expression
\be{237a}
V_{k,j} = \sum_{a_1 + \ldots + a_k = j} \left(\begin{matrix} j\\ a_1 \,\, a_2  \ldots  a_k\end{matrix}\right)    \frac{1}{(a_1+1) \cdots (a_k+1)}
\ee
where the variables $a_i$ run over non-negative integers.

\begin{thm}
\label{S2.7} {Let
$\mathcal{H}$ be an $(\alpha, c,D, \nu, \delta)$-family 
and $j \in \Bbb{N}_0$.  Then there exists a positive number
$\eta$ such that for all $h\in\cal H$ and suitable real polynomials $p_{h,j}$
of degree at most $k+j$ one has
}
 \be{quantity}
  \sum_{\substack{x_j \leq X\\1\le j \le  k}} 
\frac{(\log\langle\bfx\rangle)^j}{\langle\bfx\rangle^{\alpha}} h(\bfx) 
= p_{h,j}(\log X) + O(X^{- \eta}). 
\ee
If $c_h\neq 0$, the degree
of $p_{h,j}$ is $k+j$, and its leading coefficient is $ \alpha^k c_h V_{k, j}$.
\end{thm}

{\em Proof.} First observe that   
\begin{displaymath}
  \frac{\partial}{\partial x} \frac{\log(x y)^j}{(xy)^{\alpha}} =  \frac{-\alpha q_j(\log(x y))}{x^{\al+1}y^{\alpha}} 
\end{displaymath}
where $q_j(t) = t^j + j\al^{-1}t^{j-1}$ is a monic polynomial of degree $j$. 
Repeated use of this identity and partial summation applied to the $k$ sums 
over $\bfx$ produce the identity
\begin{equation}\label{sumxi}
  \sum_{\substack{x_j \leq X\\ 1\le j \le  k}} \frac{(\log\langle\mathbf{x}\rangle)^j 
  }{\langle\mathbf{x}\rangle^{\alpha}}h(\mathbf{x}) = 
\sum_{{\cal N}\subset\{1,\dots,k\}} \alpha^n \Xi_{\cal N}
\end{equation}
where $n=\#\cal N$, where
$$ \Xi_\emptyset =  (k\log X)^j X^{-\alpha k}\sum_{\substack{x_i \leq X\\ 1\le i\le k}} h(\bfx), $$
and when ${\cal N}$ is a nonempty subset of $\{1,\ldots,k\}$,
\be{N2}
 \Xi_{\cal N} =   X^{-\alpha m}\int_{[1, X]^n} 
\frac{q_{\cal N}(\log X^{m} \langle\bft_{\cal N}\rangle)}{\langle\bft_{\cal N}\rangle^{\alpha+1}} 
\sum_{\substack{x_i \leq t_i\\ i\in{\cal N}}} \sum_{\substack{x_i \leq X\\ i \not\in {\cal N}}} h(\mathbf{x}) 
\,\dd\mathbf{t}_{\cal N},
\ee
in which $m=k-n$, $\bft_{\cal N}=(t_i)_{i\in\cal N}$ and $q_{\cal N}$ is a certain
monic polynomial of degree $j$. 

The strategy is now to prove that for any
subset $\cal N$ of $\{1,\ldots,k\}$ there is a polynomial 
$p_{\cal N}=p_{h,\cal N}$
with the property that
\be{237} \Xi_{\cal N} =   p_{\cal N}(\log X) + O(X^{-\eta}). \ee
Once this is established, it follows from \rf{sumxi} that the 
asymptotic relation \rf{quantity} holds with {\em some} polynomial
$p_{h,j}$. The degree and leading coefficient
can then be computed by the following trick: for ${\cal K}=\{1,\dots,k\}$,
 one may use (I) in the definition of $\Xi_{\cal K}$ to deduce
that 
$$ \Xi_{\cal K} = c_h 
\int_{[1, X]^k} \frac{q_{\cal K}(\log\langle\mathbf{t}\rangle)}{\langle\mathbf{t}\rangle}\,\dd\mathbf{t} 
 + O\Bigl( \int_{[1, X]^k} \frac{(\log\langle\mathbf{t}\rangle)^j}{\langle\mathbf{t}\rangle \underset{1 \leq i \leq k}{\min} t_i^{\delta}} \,\dd\mathbf{t}\Bigr). 
$$
Since $q_{\cal K}$ is monic of degree $j$, it follows that
$$ \Xi_{\cal K} = c_h  
 \int_{[1, X]^k} \frac{(\log\langle\mathbf{t}\rangle)^j}{\langle\mathbf{t}\rangle} \,\dd\mathbf{t} 
 + O\bigl((\log X)^{k+j-1}  \bigr). 
$$
One may expand the logarithm, using the multinomial theorem. Then, an elementary
calculation yields
\begin{displaymath}\begin{split} \Xi_{\cal K}& = c_h
\sum_{a_1 + \ldots + a_k = j} \left(\begin{matrix} j\\ a_1 \,\, a_2  \ldots  a_k\end{matrix}\right)    \int_{[1, X]^k} \frac{(\log t_1)^{a_1} \cdots (\log t_k)^{a_k}} {t_1 \cdots t_k} \,\dd\mathbf{t}  + O((\log X)^{k+j-1})\\
& =   c_h  \sum_{a_1 + \ldots + a_k = j} 
\left(\begin{matrix} j\\ a_1 \,\, a_2  \ldots  a_k\end{matrix}\right)   
 \frac{1}{(a_1+1) \cdots (a_k+1)} (\log X)^{k+j}  + O((\log X)^{k+j-1}).
\end{split}\end{displaymath}
Similarly, but using (I) only as an upper bound, it is immediate that
$\Xi_{\cal N} \ll (\log X)^{n+j}$ holds for all subsets $\cal N$ of $\{1,\ldots,k\}$. 
Hence, by \rf{sumxi} and \rf{237a}, the sum on the left 
hand side of \rf{quantity} equals $\al^kV_{k,j} c_h (\log X)^{k+j}+ O((\log X)^{k+j-1})$.
This is compatible with \rf{quantity} only when $p_{h,j}$ has degree at most $k+j$,
and for $c_h\neq 0$, the degree must be $k+j$ with leading coefficient $\al^kV_{k,j} c_h$,
as required to complete the proof of
Theorem \ref{S2.7}. 

There remains the task to confirm \rf{237}. Two cases are easy. For
 ${\cal K}=\{1,\dots,k\}$, the desired expansion for $\Xi_{\cal K}$ follows from
Lemma \ref{L2.6}. Also, by (I) and the definition of $\Xi_\emptyset$,
one finds that
$$ \Xi_\emptyset = (k\log X)^j \bigl( c_h + O(X^{-\delta})\bigr) $$
confirming \rf{237}
in the case where ${\cal N}=\emptyset $.

This leaves the cases where $1\le n\le k-1$. For these $n$ we write
 $\Xi_n=\Xi_{\{1,\ldots,n\}}$.  Let $\bft=(t_1,\ldots,t_n)$.
By multinomial expansion,
there are
certain constants $\gamma_{r, s} \in \Bbb{R}$ such that
\begin{displaymath}
  q_{\{1,\ldots,n\}}(\log(\langle\bft\rangle X^{m})) = \sum_{r + s \leq j} \gamma_{r, s}
( \log\langle\bft\rangle)^r (\log X)^s,
\end{displaymath}
and consequently, one may rewrite the defining equation \rf{N2} for $\Xi_n$ as
\be{N3}
  \Xi_n =  X^{-\alpha m}\sum_{r + s \leq j} \gamma_{r, s} (\log X)^s 
\int_{[1, X]^n} \frac{ (\log \langle\mathbf{t}\rangle)^r }{\langle\mathbf{t}\rangle^{\alpha+1}} 
\sum_{\substack{x_i \leq t_i\\ 1\leq i \leq n}} \;\sum_{\substack{x_i \leq X\\n < i \leq k}} 
h(\mathbf{x}) \,\dd\mathbf{t}. 
\ee

The asymptotic evaluation of $\Xi_n$ is performed in a  manner similar to
 the proof of Lemma \ref{L2.4}. However, the details are somewhat different
because the two innermost sums in \rf{N3} are both nonempty, so that there is
at least one ``long'' sum $x_k\le X$ involved. The partition of the set
$\Delta^{(n)}(X)$ in the argument below will therefore be different from
the arrangement in the proof of Lemma \ref{L2.4}.

The first step is  to sort the $t_i$ by size. This can be done
by the argument used to prove Lemma \ref{L2.6}. Indeed, when $\sgm\in S_n$, 
we define $\sgm'\in S_k$ by $\sgm'(i)=\sgm(i)$, for $1\le i\le n$, and
$\sgm'(i)= i$ for $n < i\le k$. Then, as in the proof of Lemma \ref{L2.6},
\be{N4} \int_{[1, X]^n} \frac{ (\log \langle\mathbf{t}\rangle)^r }{\langle\mathbf{t}\rangle^{\alpha+1}} 
\sum_{\substack{x_i \leq t_i\\ 1\leq i \leq n}} \;\sum_{\substack{x_i \leq X\\n < i \leq k}} 
h(\mathbf{x}) \,\dd\mathbf{t}
=
\sum_{\sgm\in S_n} 
\int_{\Delta^{(n)}(X)} \frac{ (\log \langle\mathbf{t}\rangle)^r }{\langle\mathbf{t}\rangle^{\alpha+1}} 
\sum_{\substack{x_i \leq t_i\\ 1\leq i \leq n}}\; \sum_{\substack{x_i \leq X\\n < i \leq k}} 
h_{\sgm'}(\mathbf{x}) \,\dd\mathbf{t}.
\ee

We now construct a dissection of $\Delta^{(n)}(X)$.
Let $\beta$ be defined by \rf{defb}, and consider the intervals
$$ {\cal I}_i = (X^{\beta^{i+1}}, X^{\beta^{i}}] \quad (0\le i < n), \quad 
{\cal I}_n = [1, X^{\beta^{n}}]. $$
that provide a partition of $[1,X]$ into $n+1$ subsets. By the box principle, 
for any $\mathbf t \in \Delta^{(n)}(X)$, there is at least one ${\cal I}_i$
that contains none of the coordinates of $\mathbf t$, and the smallest such
$i$ is denoted by $i(\mathbf t)$. Once $i(\mathbf t)$ is determined, we put
$l(\mathbf t)=0$ if $t_1>X^{\beta^{i(\mathbf t)}}$, and otherwise we take $l(\mathbf t)$
to denote the largest $l$ with $t_l \le X^{\beta^{i(\mathbf t)+1}}$. Now put
$$  \Delta^{(n)}_{i,l}(X)= \{\mathbf t\in \Delta^{(n)}(X): i(\mathbf t)=i,\;
l(\mathbf t)=l\}. $$
Note that whenever $i=i(\mathbf t)$, then all the intervals ${\cal I}_0,
\ldots, {\cal I}_{i-1}$ will contain at least one coordinate of $\mathbf t$.
Hence, $ \Delta^{(n)}_{i,l}(X)$ will be nonempty only when $i+l\le n$,
and $\Delta^{(n)}(X)$ is the disjoint union of these sets.
We now write 
\be{defJ}
  J_{r,n ,i,l}(h) =   X^{-\alpha m} \int_{\Delta^{(n)}_{i,l}(X)} 
\frac{ (\log\langle\mathbf{t}\rangle)^r }{\langle\mathbf{t}\rangle^{\alpha+1}}  
\sum_{\substack{x_j \leq t_j\\ 1\leq j \leq n}} \;
\sum_{\substack{x_j \leq X\\n < j \leq k}} h(\mathbf{x}) \,\dd\mathbf{t}. 
\ee
Then, by \rf{N3} and \rf{N4},
\be{sumJ}
  \Xi_n =  \sum_{\sgm\in S_n} \sum_{r + s \leq j} \gamma_{r, s} (\log X)^s  
\sum_{i+l\le n} J_{r,n,i ,l}(h_\sgm).
\ee  

We begin with the  evaluation of $J_{r,n,i ,0}(h)$. Note that $\Delta^{(n)}_{0,0}(X)$
is empty so that we may suppose that $i\ge 1$. Then, for $\mathbf t \in
\Delta^{(n)}_{i,0}(X)$, one has $t_1\ge X^{\beta^{i}} \ge X^{\beta^{n}}$, and (I)
delivers
$$  \sum_{\substack{x_j \leq t_j\\ 1\leq j \leq n}} \;
\sum_{\substack{x_j \leq X\\n < j \leq k}} h(\mathbf{x}) =
c_h X^{\alpha m} \langle\mathbf{t}\rangle^{\alpha} + 
O(X^{\alpha m -\del\beta^n} \langle\mathbf{t}\rangle^{\alpha}).$$
By \rf{defJ} and straightforward estimates,
\be{Jint}  J_{r,n,i ,0}(h) =  c_h 
\int_{\Delta^{(n)}_{i,0}(X)} \frac{ (\log\langle\mathbf{t}\rangle)^r }{\langle\mathbf{t}\rangle }  \,\dd \mathbf{t} + O(X^{ -\delta\beta^n/2 } ). \ee

It remains to evaluate the integral on the right hand side here. We claim that
there is a polynomial $T$ depending only on $n,i,\beta$ and $r$ such that
\be{poly0}\int_{\Delta^{(n)}_{i,0}(X)} \frac{ (\log\langle\mathbf{t}\rangle)^r }{\langle\mathbf{t}\rangle }  
\,\dd \mathbf{t} = T(\log X).
\ee
To see this,
let
$n=u_1>u_2>\ldots>u_i>u_{i+1}=1$ be a collection of natural numbers, and
put $\bfu=(u_1,\ldots, u_{i+1})$. Let
$$\Gamma_{i,\bfu}=\{\mathbf t\in\Delta^{(n)}(X) : t_\rho\in{\cal I}_\lambda \text{ for }
u_{\lambda+2}< \rho \le u_{\lambda+1}\; (0\le\lambda <i)\}.$$  
By construction, $\Delta^{(n)}_{i,0}(X)$ is the disjoint union of the
$\Gamma_{i,\bfu}$. Hence, by multinomial expansion, the integral in \rf{poly0}
equals
$$ \sum_{n=u_1>u_2>\ldots>u_i>u_{i+1}=1}\,
\sum_{a_1 + \ldots + a_n = r} \left(\begin{matrix} j\\ a_1 \,\, a_2  \ldots  a_k\end{matrix}\right)
\int_{\Gamma_{i,\bfu}} \frac{(\log t_1)^{a_1} \cdots (\log t_n)^{a_n}} {t_1 \cdots t_n} \,\dd\mathbf{t}.$$
By definition of $\Gamma_{i,\bfu}$, this last integral factorises into integrals
over $t_\rho\in{\cal I}_\lambda$ with $u_{\lambda+2}< \rho \le u_{\lambda+1}$
of dimension $u_{\lambda+1}-u_{\lambda+2}$. A typical such integral takes the shape
$$ \int_{Y^{\beta}\le v_1<\ldots <v_s\le Y} \frac{(\log v_1)^{b_1} \cdots (\log v_s)^{b_s}} {v_1 \cdots v_s} \,\dd\mathbf{v} $$
where $s$ is the dimension, the $b_\rho$ are some of the $a_1,\ldots, a_n$,
and $Y=X^{\beta^{\lambda}}$ for some $\lambda < n$.
This integral can be computed explicitly, and is then seen to be a polynomial
in $\log Y$, and hence also a polynomial in $\log X$. Collecting together confirms \rf{poly0}.  

By \rf{Jint} and \rf{poly0}, we see that $J_{r,n,i ,0}(h)$ is a polynomial
in $\log X$, up to an acceptable error $O(X^{-\delta\beta^n/2})$. By (III), the same is true for
$J_{r,n,i ,0}(h_\sgm)$. 

\smallskip

The next case we consider is $l=n$. Since $l+i\le n$, this forces $i=0$, and one 
readily checks from the relevant definitions that
$\Delta^{(n)}_{0,n}(X) =  \Delta^{(n)}(X^{\beta})$. Recall that $\beta \le \nu$. Hence, on writing
$\bfx'=(x_1,\ldots,x_l)$,  we deduce from (II) that
$$ \sum_{\substack{x_j \leq X\\n < j \leq k}} h(\mathbf{x}) =
c_{h,l}(\bfx')X^{\al m}+ O(X^{\al m-\del} |\bfx'|^D). $$
For $\mathbf t\in\Delta^{(n)}(X^{\beta})$ we may sum over $\bfx'\le \mathbf t'$ to infer that
$$  \sum_{\substack{x_j \leq t_j\\ 1\leq j \leq n}} \;
\sum_{\substack{x_j \leq X\\n < j \leq k}} h(\mathbf{x}) =
X^{\al m} \sum_{\bfx'\le \mathbf t} c_{h,l}(\bfx') + O(X^{\al m-\del}|\mathbf t|^{n+D}). $$
This may be injected into \rf{defJ}. Then, recalling that $\beta\le \del/(4n+2D)$, one concludes
that
$$ J_{r,n,0,n}(h) =  \int_{\Delta^{(n)}(X^{\beta})} 
\frac{ (\log\langle\mathbf{t}\rangle)^r }{\langle\mathbf{t}\rangle^{\alpha+1}}
\sum_{\bfx'\le \mathbf t} c_{h,l}(\bfx')\,\dd\mathbf t + O(X^{-\del/2}). $$
By Lemma \ref{L2.2}, we may apply Lemma \ref{L2.4} with $c_{h,l}$ in place of $h$ to
evaluate the integral on the right hand side here. It follows that $ J_{r,n,0,n}(h) $
is equal to a polynomial in $\log X$, up to an error not excceding
$O(X^{-\eta})$, for some suitable $\eta>0$. 

\smallskip

We are left with the cases where $1\le l< n$. Write 
$\mathbf{t} = (\mathbf{t}', \mathbf{t}'')$ 
with $\mathbf{t}' = (t_1, \ldots, t_l)$ 
and $\mathbf{t}'' = (t_{l+1}, \ldots, t_n)$.
Another inspection of the relevant definitions shows that $\mathbf t \in\Delta^{(n)}_{i,l}(X)$
holds if and only if  $\mathbf t' \in\Delta^{(l)}(X^{\beta^{i+1}})$ and 
$\mathbf t'' \in\Delta^{(n-l)}_{i,0}(X)$. Moreover, since $\beta\le\nu$, one may use (II)
to confirm that whenever $\mathbf t \in\Delta^{(n)}_{i,l}(X)$, then
$$ \sum_{\substack{x_j \leq t_j\\ l< j \leq n}} \;
\sum_{\substack{x_j \leq X\\n < j \leq k}} h(\mathbf{x}) = c_{h,l}(\bfx')\langle \mathbf t''
\rangle^\al X^{\al m} + O(X^{\al m} \langle \mathbf t''\rangle^\al t_{l+1}^{-\del} |\mathbf x'|^D). $$
We sum over $\bfx'\le \mathbf t'$ and recall that $t_{l+1}\ge X^{\beta^{i}}$. This produces
$$ \sum_{\substack{x_j \leq t_j\\ 1\leq j \leq n}} \;
\sum_{\substack{x_j \leq X\\n < j \leq k}} h(\mathbf{x}) =
X^{\al m}\langle \mathbf t''\rangle^\al \sum_{\bfx'\le \mathbf t'} c_{h,l}(\bfx')
+ O(X^{\al m-\del \beta^{i}} \langle \mathbf t''\rangle^\al |\mathbf t'|^{n+D}).
$$
Now multiply with $ \langle \mathbf t\rangle^{-\al-1} (\log  \langle \mathbf t\rangle)^r$
and integrate over $\Delta^{(n)}_{i,l}(X)$. The error term above then becomes
\begin{displaymath}\begin{split}
X^{\al m-\del \beta^{i}} &(\log X)^r \int_{\Delta^{(n)}_{i,l}(X)} \langle \mathbf t\rangle^{-1} |\mathbf t'|^{n+D} \,\dd\mathbf t \\ & \mbox{} \ll
 X^{\al m-\del \beta^{i}} (\log X)^{r+n} \int_{\Delta^{(l)}(X^{\beta^{i+1}})} |\mathbf t'|^{n+D} \,\dd\mathbf t'
\ll X^{\al m-\del\beta^n/2}.
\end{split}\end{displaymath}
In the last inequality, we have used again that $\beta\le \del/(4n+2D)$. We insert the results into
\rf{defJ} and apply the binomial theorem to conclude that, up to an error $ O(X^{-\del\beta^n/2})$, 
the expression $J_{r,n ,i,l}(h)$ equals
$$    \sum_{r'+r''=r}\left(\begin{matrix} r\\r' \end{matrix}\right) 
\int_{\Delta^{(n-l)}_{i,0}(X)} 
\frac{(\log \langle \mathbf t''\rangle)^{r''}}{\langle \mathbf t''\rangle}
\,\dd\mathbf t'' 
\int_{\Delta^{(l)}(X^{\beta^{i+1}})}  
\frac{(\log \langle \mathbf t'\rangle)^{r'}}{\langle \mathbf t''\rangle^{\al+1}}
\sum_{\bfx' \le \mathbf t'} c_{h,l}(\bfx') \,\dd\mathbf t' . $$
By \rf{poly0}, the first integral is a polynomial in $\log X$.
The second integral coincides with a suitable polynomial in $\log X$,
up to an error not exceeding $O(X^{-\eta})$ for some suitable number $\eta>0$.
This follows once again from  Lemma \ref{L2.2} and Lemma \ref{L2.4}.

To sum up the above deliberations, we have now shown that $J_{r,n ,i,l}(h)$
equals a suitable polynomial in $\log X$,
up to an error not excceding $O(X^{-\eta})$ for some suitable number $\eta>0$,
for all revelant parameters $r,n ,i,l$. By (III), the same is true for  $J_{r,n ,i,l}(h_\sgm)$,
and \rf{237} for ${\cal N}=\{1,\ldots,n\}$ now follows from \rf{sumJ}. By (III) again,
this confirms \rf{237} for all $\cal N$. The proof of
Theorem \ref{S2.7} is complete.

\subsection{Away from the spikes}
In this section, we discuss the contribution to the
sum \rf{211} where all the variables $u_j$ are rather large.
More precisely, we choose a threshold $\bfW\in[1,N]^k$ and consider
\be{251}  \Ups(N,\bfW) = 
\multsum{w_1 w_2\cdots w_k\le N}{\bfw > \bfW} h(\bfw).
\ee
When all entries of $\bfW$ are reasonably large,
the variables of summation stay away from the spiky part of
the hyperbolic constraint $\langle\bfw\rangle\le N$, so that 
one would hope to handle 
this sum based on the condition (I) alone. This is indeed the case. 
The asymptotic formula features the polynomial
\be{255} p_k(t) = \sum_{l=0}^{k-1} \f{(-1)^{k+1+l}}{l!} t^l. \ee

\begin{lem}
\label{L2.9} {  
Let $\cal H$ be a set of functions $h:\NN^k \to [0,\infty)$
that satisfies the condition {\rm (I)} with respect to $(\al,c,\del)$, 
and suppose that the threshold
satisfies $ \langle\bfW\rangle\le N^{1/2}$ and $\min W_j \ge (\log N)^{2k/\del}$. 
Then}
$$  \Ups(N,\bfW) = c_h N^\al p_k\Big(\al\log\f{ N}{\langle \bfW\rangle}\Big) + O(N^{\al}(\min W_j)^{-\del/(2k)} (\log N)^k).$$
\end{lem}

The proof depends on the following combinatorial
identity.

\begin{lem}
\label{L2.10} {Let $k$ and $J$ be natural numbers. Then 
for $t\in\CC$ one has
\be{252} (1-t)^k
\multsum{j_1+\ldots+j_k\le J}{j_r\ge 0} t^{j_1+\ldots+ j_k} =
1 - t^{J+1} \sum_{l=0}^{k-1}  \Big({ J+l \atop l} \Big)(1-t)^{l} . 
\ee}
\end{lem}

When  $k=1$, the claim in Lemma \ref{L2.10} is
 the
familiar evaluation of the geometric sum. We proceed by induction and 
suppose that the formula is 
known for $k-1$ in place of $k$.  Then the left hand side of \rf{252}
equals
\begin{eqnarray*}
\lefteqn{ 
(1-t)^k \sum_{j=0}^J t^j \sum_{j_1+\ldots+j_{k-1}\le J-j} t^{j_1+\ldots+ j_{k-1}}
} \\
&=& (1-t) \sum_{j=0}^J t^j \Big( 1 - t^{J+1-j} 
\sum_{l=0}^{k-2} \Big({ J+l-j \atop l} \Big)(1-t)^l \Big) .
\end{eqnarray*}
Now replace $l+1$ by $l$ in the inner sum. The above then becomes
\be{253}
 1-t^{J+1} - t^{J+1} \sum_{j=0}^J \sum_{l=1}^{k-1}
 \Big({J+l-1-j \atop l-1}\Big) (1-t)^l .
\ee
However,
$$ \sum_{j=0}^J \Big({J+l-1-j \atop l-1}\Big)
=\sum_{j=0}^J \Big({l-1+j \atop l-1}\Big) = \Big({J+l \atop l}\Big),
$$
as one may verify by induction on $J$. Now \rf{252} follows from
\rf{253}.

\medskip

{\em Proof of Lemma} \ref{L2.9}. We begin with an enveloping argument
to reduce the evaluation of $\Ups(N,\bfW)$ to box sums of the
type
\be{256} H(\bfU^+,\bfU^-) = \sum_{\bfU^- < \bfu \le \bfU^+} h(\bfu).\ee
Let $\Th$ be a real number, $J$ be a natural number, and
suppose that $1<\Th<3$ and $\Th^J = N/ \langle\bfW\rangle$. 
We shall optimize $J$
later, but already note that permissible values of $J$ satisfy $J
\gg \log N$.  For $j\ge 0$, let $U_{r,j}=W_r\Th^j$, and define
$$ \bfU_\bfj = (U_{1,j_1},\ldots,U_{k,j_k}). $$
Let ${\bf 1}= (1,1,\ldots,1)$. We consider boxes $\bfU_\bfj < \bfu \le \bfU_{\bfj +\bfone}$. This box
lies inside the range of summation $u_1\cdots u_k\le N$ in \rf{251}
whenever
$ U_{1,j_1+1}\cdots U_{k,j_k +1} \le N$ which in turn holds if and only
if $|\bfj|_1 \le J-k$. Here, and later in this proof, we write
$$ |\bfj|_1 = j_1+ \ldots + j_k. $$
In the opposite direction, let $\bfu$ be a point with $\bfu > \bfW$
and $u_1\cdots u_k\le N$. Then, there is a unique $\bfj\in\NN_0^k$
with $\bfU_\bfj < \bfu \le \bfU_{\bfj+\bfone}$. The inequalities 
$$  U_{1,j_1}\cdots U_{k,j_k} < u_1\cdots u_k \le N$$
imply  $|\bfj|_1\le J$. By \rf{251} and \rf{256}, these considerations
show that
\be{257} 
\sum_{|\bfj|_1\le J-k} H(\bfU_{\bfj+\bfone},\bfU_\bfj) \le \Ups(N,\bfW)
\le \sum_{|\bfj|_1\le J} H(\bfU_{\bfj+\bfone},\bfU_\bfj).
\ee
By the inclusion-exclusion principle, the sum $H(\bfU^+,\bfU^-)$ can be
expressed in terms of the sum $H(\bfX)$ that was defined in \rf{212}.
In the special case needed here, this strategy gives the identity
$$  H(\bfU_{\bfj+\bfone},\bfU_\bfj) = \sum_{\bfs\in\{0,1\}^k} (-1)^{k-|\bfs|_1}
H(\bfU_{\bfj+\bfs})
$$
whence by (I), it now follows that
$$ H(\bfU_{\bfj+\bfs}) = c_h \langle\bfW\rangle^{\al} \Th^{\al|\bfj+\bfs|_1}
+O(  \langle\bfW\rangle^{\al}  (\min  W_j)^{-\del} \Th^{\al|\bfj|_1}). $$
 Now put 
$$  J^+=J , \quad J^- = J-k $$
and study the sums 
$$ \Ups^\pm =\sum_{|\bfj|_1\le J^\pm} H(\bfU_{\bfj+\bfone},\bfU_\bfj)$$
that occur in the sandwich inequalities \rf{257}. Combining the
preceding displays yields
$$ \Ups^\pm = c_h \langle\bfW\rangle^{\al} 
\sum_{|\bfj|_1\le J^\pm}
 \sum_{\bfs\in\{0,1\}^k} (-1)^{k-|\bfs|_1} \Th^{\al|\bfj+\bfs|_1} 
+  O\Big(  \langle\bfW\rangle^{\al}  (\min  W_j)^{-\del} 
\sum_{|\bfj|_1\le J^\pm}\Th^{\al|\bfj|_1} \Big) . $$
The obvious identity
$$ \sum_{\bfs\in\{0,1\}^k} (-1)^{k-|\bfs|_1} T^{|\bfs|_1} = (T-1)^k $$
allows us to rewrite this as
$$ \Ups^\pm = c_h (\Th^\al -1)^k \langle\bfW\rangle^{\al} \sum_{|\bfj|_1\le J^\pm}
\Th^{\al|\bfj|_1} + O\Big(  \langle\bfW\rangle^{\al}  (\min  W_j)^{-\del} 
\sum_{|\bfj|_1\le J^\pm}\Th^{\al|\bfj|_1} \Big). $$
It will now be convenient to define $r$ by $\min W_j = W_r$. Then, 
the error term above does not exceed
$$  \langle\bfW\rangle^{\al} W_r^{-\del} \Th^{\al J} 
\#\{\bfj: |\bfj|_1\le J\} \ll N^{\al} W_r^{-\del}J^k.$$
To compute the leading term, we multiply \rf{252} with $(-1)^k$ and choose 
$t=\Th^\al$. Then $t^{J^+} = N^{\al} \langle\bfW\rangle^{-\al}$, and one finds that
\be{258} \Ups^+ = c_h N^\al\Th^\al \sum_{l=0}^{k-1} \Big({J+l \atop l}\Big)
(-1)^{k+1+l} (\Th^\al-1)^l + O( \langle\bfW\rangle^{\al}+ N^{\al}W_r^{-\del}J^k). \ee
The same argument also gives
\be{259} \Ups^- = 
c_h N^\al\Th^{\al(1-k)} \sum_{l=0}^{k-1} \Big({J-k+l \atop l}\Big)
(-1)^{k+1+l} (\Th^\al-1)^l + 
O( \langle\bfW\rangle^{\al}+ N^{\al}W_r^{-\del}J^k). \ee

We now choose $J=[W_r^{\del/(2k)}]$ so that $J\ge \log N$ and  $1<\Th\le\mathrm e$,
as required. We then have
$$ \Th = \exp\Big(J^{-1} \log \f{N}{\langle\bfW\rangle}\Big) = 
1 + J^{-1} \log \f{N}{\langle\bfW\rangle}
+ O(J^{-2}(\log N)^2). $$
Binomial expansion gives
\be{2511}
\Th^\al = 1 + \al J^{-1}  \log \f{N}{\langle\bfW\rangle}
+ O(J^{-2}(\log N)^2). \ee
One also has
$$ \Big( {J+l \atop l} \Big) = \f{J^l}{l!} + O(J^{l-1}) $$
so that
$$ \Big({J+l \atop l}\Big)(\Th^\al-1)^l = \f{\al^l}{l!} 
\Big(\log \f{N}{\langle\bfW\rangle}\Big)^l
+ O(J^{-1}(\log N)^{l+1}). $$
We multiply with $(-1)^{k+l+1}$ and sum over $l$. Recalling the notation
introduced in \rf{255}, we then infer from \rf{258}, \rf{2511} and the 
preceding display that
\be{2512}
\Ups^+ = c_h N^\al p_k\Big(\al\log \f{N}{\langle\bfW\rangle}\Big) +
O(N^\al(\log N)^k J^{-1} + \langle\bfW\rangle ^{\al} + N^{\al}W_r^{-\del} J^k).
\ee
A simple cosmetical change in this argument, now starting from
\rf{259} provides the
same asymptotic formula for  $\Ups^-$. Consequently, by \rf{257},
this formula also holds for $\Ups(N,W)$. Our choice for $J$ then yields
the assertion of Lemma \ref{L2.9}.

\subsection{The endgame} We are ready to
assemble the puzzle. The contribution to the sum \rf{211} resulting from 
summands $h(\bfx)$ where all coordinates $x_j$ are ``large'' can be evaluated 
by Lemma \ref{L2.9}. The contribution from terms where all $x_j$ are ``small'' 
will not be of significance, for obvious reasons. This leaves summands where
$\bfx$ has small and large coordinates simultaneously. Here the strategy is to
sum over the large coordinates with the aid of Lemma \ref{L2.9}, and then sum
the result over the small coordinates by an appeal to Theorem \ref{S2.7}. This is
successful only if there is a huge gap between the small and the large
coordinates. One can always find such a gap, but its position will depend
on $\bfx$. Hence, we will follow a strategy that is largely similar to 
the one used in the proof of Theorem \ref{S2.7}, but there are additional 
combinatorial complications because summands which have at least two
equal coordinates in the vector $\bfx$ affect the lower order terms
in $P_h$.   

Fix an $(\al,c,D,\nu,\delta)$-family $\cal H$ and $k\ge 2$. The set $\cal H$ 
remains such a family if the values 
of $\nu$ and $\delta$ shrink, and we may therefore suppose that
$$ \nu\le 1/(2k), \quad \delta \le 1/2.  $$
Also, as in Lemma \ref{L2.3}, put $A=D+(k+1)\al + \nu^{-1}(1+\al)$ and 
then define
\be{defB} B= 4Ak^2 / \delta. \ee
Fix a parameter $V\ge 4$, and suppose that $N$ is sufficiently large for
 $V^{B^k}\le N^\nu$ to hold.  Then use the sequence
\be{defV} V_0=1, \quad V_1=V, \quad V_j = V_{j-1}^B \;\;(2\le j\le k), \quad V_{k+1}=N \ee
to define the intervals
$$ {\cal V}_0= [V_0, V_1],\quad {\cal V}_j = (V_j, V_{j+1}] \quad (1\le j\le k). 
$$
Then $[1,N)$ is the disjoint union of the $k+1$ sets ${\cal V}_j$. Hence, for
$\bfx\in\NN^k$ with $\langle\bfx\rangle \le N$ there is at least one 
$l\in\{0,\ldots,k\}$ with $x_j\not\in {\cal V}_l$ for all $1\le j\le k$, 
and we may define
$l(\bfx)$ as the largest such $l$.  
We sort terms in \rf{211} according to the value of $l(\bfx)$. Thus, we write
\be{261} \Ups_l(N) = \multsum{\langle \bfx\rangle \le N}{l(\bfx)=l} h(\bfx)
\ee
and note that
\be{262} \Ups(N) = \Ups_0(N) + \Ups_1(N) +\ldots +\Ups_k(N). \ee
The condition that $l(\bfx) =k$ is equivalent with $|\bfx|\le V_k$, whence by 
(I)  and \rf{defB} one obtains the crude bound
\be{263} \Ups_k(N)\le H(V_k,\ldots,V_k) \ll V_k^{k\al} \ll N^{\al/2}. \ee

In line with the comments preceding the current discussion, in our approach to 
the estimation of $\Ups_l(N)$ we will treat coordinates $x_j$ with $x_j\le V_l$
as ``small'', and all $x_j>V_{l+1}$ as ``large''. To make this precise, let 
$0\le l\le k-1$ and $\bfx\in\NN^k$ with $\langle \bfx\rangle \le N$ and 
$l(\bfx)=l$. With such an $\bfx$, we associate the sets
\be{264} {\cal J}=\{j: \, 1\le j\le k,\, x_j\le V_l\}, \quad
{\cal L}_m = \{j:\,x_j\in{\cal V}_m\}.\ee
Then, the maximality of $l(\bfx)$ implies that ${\cal L}_{l+1},\ldots,
{\cal L}_k$ are non-empty, and $\{1,\ldots,k\}$ is the disjoint union of
${\cal J}$ and ${\cal L}_m$, $l<m\le k$.   

An ensemble of sets ${\fk A}=\{ {\cal J}, {\cal L}_{l+1}, {\cal L}_{l+2},\ldots,
 {\cal L}_{k}\}$ with all $ {\cal L}_{m}$ non-empty and such that $\{1,\ldots,k\}$
is the disjoint union of ${\cal J}, {\cal L}_i$ $(l<i\le k)$    will be referred to as {\em permissible
to the natural number $l$}. For any $\bfx\in\NN^k$ with $\langle \bfx\rangle \le N$ and 
$l(\bfx)=l$ the sets \rf{264} form an ensemble permissible to $l$ that we denote by
${\fk A}(\bfx)$. For an ensemble $\fk A$ that is permissible to some $l\in\{0,\ldots,k-1\}$,
we now define
$$ \Ups_{\fk A}(N) = \multsum{\langle \bfx\rangle \le N}{{\fk A}(\bfx)=\fk A} h(\bfx) 
$$
and then have
\be{266} \Ups_l(N) = \sum_{\fk A}  \Ups_{\fk A}(N)\ee
where the sum extends over all ensembles $\fk A$ that are permissible to $l$.

With \rf{262} and \rf{266} in hand, our initial decomposition of $\Ups(N)$ is now complete.
We proceed to deduce Theorem \ref{S2.1} in two steps, each depending on this decomposion yet with a
different choice for the parameter $V$. We begin with a weak form of Theorem \ref{S2.1} in which
a leading term is already identified.

\begin{lem}
\label{L2.11} { Let $\mathcal{H}$ be an $(\alpha, c,D, \nu, \delta)$-family. 
Then, uniformly for any $h\in\cal H$, }
\be{266a} \Ups(N) = \f{c_h \al^{k-1}}{(k-1)!} N^\al (\log N)^{k-1} + O(N^\al
(\log N)^{k-2} \log\log N). \ee
\end{lem}

{\em Proof}. The case $k=1$ follows from (I). Hence, we may assume that $k\ge 2$. We choose
\be{267}  V= (\log N)^{B} \ee
with $B$ as in \rf{defB} and decompose $\Ups(N)$ according to \rf{262} and \rf{266}.
By \rf{263}, the summand $\Ups_k(N)$ is absorbed into the error term in \rf{266a}. Next,
consider the ensemble ${\fk A} = \{\emptyset, \{1,2,\ldots,k\}\}$ that is permissible to 
$k-1$, and note that $\Ups_{\fk A}(N) = \Ups(N,(V_k,\ldots,V_k))$, in the notation
introduced in \rf{251}. We have $V_k = (\log N)^{B^k}$ so that Lemma \ref{L2.9} yields
\begin{eqnarray}
\Ups(N,(V_k,\ldots,V_k)) & = &
c_h N^\al p_k(\al \log NV_k^{-k}) + O(N^\al) \nonumber \\
&=& c_h N^\al \f{\al^{k-1}}{(k-1)!} (\log N)^{k-1} + O(N^\al(\log N)^{k-2}\log\log N).
\label{268}
\end{eqnarray}
Here the final line corresponds to the right hand side of \rf{266a}. To complete the proof of
Lemma \ref{L2.11}, it remains to show that for any other ensemble $\fk A$ that is permissible 
to some $l$ the sum $\Ups_{\fk A}(N)$ can be absorbed into the error term in \rf{266a}.

First consider ensembles in which the set $\cal J$ is empty. Thus suppose
that ${\fk A}= \{\emptyset, {\cal L}_{l+1}, \ldots, {\cal L}_k\}$ is permissible
to $l$. The case $l=k$ is discussed in \rf{268} so that we are reduced to
the range $0\le l\le k-2$. Since ${\cal L}_{l+1}$ is non-empty, at least one
of the $x_j$ is contrained to ${\cal V}_{l+1}$, and 
we temporarily suppose that this is so for $x_1$. Any other coordinate $x_j$
of an $\bfx$ with ${\fk A}(\bfx)=\fk A$ must obey the 
inequality $x_j\ge V_{l+1}$, and consequently, one has the crude upper bound
\begin{eqnarray*}
\Ups_{\fk A} (N) &\le & \sum_{V_{l+1}<x_1 \le V_{l+2}} 
\trisum{x_j > V_{l+1}}{2\le j\le k}{\langle\bfx\rangle\le N}
h(\bfx) \\
&=& \Ups(N,(V_{l+1},V_{l+1},\ldots,V_{l+1})) -
\Ups(N,(V_{l+2},V_{l+1},\ldots,V_{l+1})) .
\end{eqnarray*}
Recall that $0\le l\le k-2$ so that $V_1\le V_{l+1} < V_{l+2} \le V_k$. Lemma \ref{L2.9}
now delivers the bound
$$ \Ups_{\fk A} (N) \le c_h N^{\al} \big( p_k(\al \log NV_{l+1}^{-k})
-p_k(\al \log NV_{l+2}^{-1}V_{l+1}^{1-k})\big) + O(N^\al) $$
which readily implies the desired estimate
$$  \Ups_{\fk A} (N) 
\ll N^\al \Big(\log\f{V_{l+2}}{V_{l+1}} \Big) (\log N)^{k-2}
\ll N^\al  (\log N)^{k-2} \log \log N. $$
By (III), this upper bound remains valid if another index takes the special
role of $j=1$ in the above argument. This completes the discussion of ensembles
with ${\cal J}=\emptyset$.

Next, suppose that ${\fk A} = \{ {\cal J},{\cal L}_{l+1},\ldots,{\cal L}_k\}$
is an ensemble permissible to $l$ with ${\cal J}$ non-empty. Since all 
${\cal L}_m$ are also non-empty, this enforces that $l\ge 1$. In view of (III),
we may suppose that ${\cal J} = \{1,\ldots,r\}$ with some $r\ge 1$. We write
$\bfw=(x_1,\ldots,x_r)$ and $\bfx=(\bfw,\bfy)$ with $y_j=x_{r+j}$. For $\bfx$ with ${\fk A}(\bfx)=\fk A$, one has $y_j > V_{l+1}$ for 
$1\le j\le k-r$ and $x_j\le V_l$ for $1\le j \le r$. It follows that
$$ \Ups_{\fk A}(N) \le 
\sum_{|\bfw|\le V_l} 
\trisum{\langle\bfy\rangle\le N/\langle\bfw\rangle}{y_j>V_{l+1}}{1\le j\le k-r}
h(\bfw,\bfy).
$$
We evaluate the inner sum  by applying Lemma \ref{L2.9} to the function $g_{h,\bfw}$
considered in Lemma \ref{L2.3}. This yields
\begin{eqnarray*}
 \trisum{\langle\bfy\rangle\le N/\langle\bfw\rangle}{y_j>V_{l+1}}{1\le j\le k-r}
h(\bfw,\bfy) &=& 
\f{c_{h,r}(\bfw)}{\langle\bfw\rangle^\al} N^\al p_{k-r} \Big(\al \log \f{NV_{l+1}^{r-k}}{\langle\bfw\rangle}\Big)
+ O(N^\al V_{l+1}^{-\delta/2k} (\log N)^k \langle\bfw\rangle^A) 
\\ &\ll&
\f{c_{h,r}(\bfw)}{\langle\bfw\rangle^\al} N^\al (\log N)^{k-r-1} +
N^\al V_{l+1}^{-\delta/2k} (\log N)^k \langle\bfw\rangle^A.
\end{eqnarray*}
One may now sum over $\bfw$ with the aid of Theorem \ref{S2.7} to deduce that
$$
\Ups_{\fk A} \ll N^\al (\log N)^{k-r-1} (\log V_l)^r +
N^\al V_l^{r(A+1)}V_{l+1}^{-\delta/2k} (\log N)^k . $$
The first term is acceptable by \rf{267}, and the second is $O(N^\al)$
in view of \rf{defB}. This completes the proof of Lemma \ref{L2.11}.

\medskip

In the proof of Lemma \ref{L2.11} it was possible to estimate $\Ups_{\fk A}(N)$
rather crudely, once the leading term was identified in \rf{268}. We now 
build up a related argument to establish Theorem \ref{S2.1}. Let $B$ be as defined 
in \rf{defB} and put $\kp = \nu/(kB^k)$. 
Now take $V=N^\kp$ to obtain another decomposition of $\Ups(N)$ via \rf{261}
and \rf{262}. We shall show that for an ensemble $\fk A$ permissible to
some $l\in\{0,\ldots,k-1\}$ there is a real polynomial $P_{\fk A} =
P_{{\fk A},h}$ of degree at most
$k-1$ and such that
\be{269} \Ups_{\fk A}(N) = N^\al P_{\fk A}(\log N) + O(N^{\al-\eta}) \ee
holds with some suitably small $\eta>0$. 
Once this is established, we deduce from \rf{261} and \rf{266} in conjunction
with \rf{263} that the asymptotic relation
$$ \Ups (N) = N^\al P^*_h (\log N) + O(N^{\al-\eta}) $$
holds with
$$ P^*_h = \sum P_{{\fk A},h} $$
in which the sum extends over all ensembles permissible to some 
$l\in\{0,\ldots,k-1\}$. Then $P^*_h$ is a polynomial of degree at most $k-1$,
but Lemma \ref{L2.11} shows that for $c_h> 0$ the degree is indeed $k-1$, and that
the leading coefficient is as claimed in Theorem \ref{S2.1}. 
Further, when $c_h=0$, then Lemma \ref{L2.11} implies that the degree of
$P^*_h$ does not exceed $k-2$.
Consequently, the verification
of \rf{269} will complete the proof of Theorem \ref{S2.1}.

\smallskip

First consider an ensemble  
${\fk A}= \{\emptyset, {\cal L}_{l+1}, \ldots, {\cal L}_k\}$ that is permissible
to $l$. If $l=k-1$ then ${\cal L_k}= \{1,\ldots,k\}$ and
$\Ups_{\fk A}(N) = \Ups(N,(V_k,\ldots,V_k))$, as observed in the proof of Lemma \ref{L2.11}.
In contrast with the discussion in \rf{268}, now $V_k$ is a fixed power of $N$,
and Lemma \ref{L2.9} delivers the asymptotic relation
$$ \Ups_{\fk A} (N) = c_h N^\al p_k(\al \log NV_k^{-k}) + O(N^{\al+\eps}V_k^{-\delta/2k}), $$
as is required in \rf{269}. It remains to consider the
case where $0\le l\le k-2$. Then the conditions on $\bfx$ in the sum defining 
$\Ups_{\fk A}(N)$ are $V_m<x_j\le V_{m+1}$ for $j\in{\cal L}_m$. For $l<m<k$ we 
use $(V_m,V_ {m+1}] = (V_m,\infty)\setminus (V_{m+1},\infty)$ and the 
inclusion-exclusion principle to obtain a representation
\be{2610}
 \Ups_{\fk A} (N) = \sum_\bfV (-1)^{\eps(\bfV)} \Ups(N,\bfV)
\ee
where $ \Ups(N,\bfV) $ is given by \rf{251}, where $\bfV$ runs through the vectors
$\bfV=(V^{(1)},\ldots, V^{(k)})$ with $V^{(j)}$ either $V_m$ or $V_{m+1}$ when $j\in{\cal L}_m$
with $m<k$,
and $V^{(j)} = V_k$ for $j\in{\cal L}_k$, and where $\eps(\bfV)\in\{0,1\}$ is chosen 
appropriately. Since $V_1=N^\kp$, it follows from Lemma \ref{L2.9} that
$$ \Ups(N,\bfV) = c_h N^\al p_k(\al \log N/\langle\bfV\rangle ) + O(N^{\al-\eta}) $$
holds for all $\bfV$ under consideration. Since all $V_j$ with $j\ge 1$ are positive 
powers of $N$, one notes that $\al\log  N/\langle\bfV\rangle$ is a constant multiple of
$\log N$. Thus, any  $ \Ups(N,\bfV) $ in \rf{2610} satisfies an asymptotic formula of 
the type desired in \rf{2610}, and so does $ \Ups_{\fk A} (N)$. This confirms \rf{269}
for ensembles in which $\cal J$ is the empty set.

\smallskip
Now let ${\fk A}= \{{\cal J}, {\cal L}_{l+1}, \ldots, {\cal L}_k\}$ be an ensemble
permissible to $l\in\{0,\ldots,k-1\}$, and suppose that ${\cal J}=\{1,\ldots,r\}$
with $r\ge 1$. As in the proof of Lemma \ref{L2.11} put $\bfx = (\bfw,\bfy)$ with 
$y_j=x_{r+j}$. Then  
$$ \Ups_{\fk A}(N) = 
\sum_{|\bfw|\le V_l} 
\sum_{\bfy}
h(\bfw,\bfy)
$$
where the sum over $\bfy$ is constrained by
$\langle\bfy\rangle \le N/\langle\bfw\rangle$ and $y_j\in{\cal V}_m$ with some
appropriate $m=m(j)\in\{l+1,\ldots,k\}$. As in the previous argument, for $m(j)<k$, 
we resolve
the condition $y_j\in{\cal V}_m = (V_m,\infty)\setminus (V_{m+1},\infty)$ by the inclusion-exclusion principle, and rewrite the previous expression as
$$ 
 \Ups_{\fk A} (N) =\sum_{|\bfw|\le V_l} \sum_\bfV (-1)^{\eps(\bfV)}
\multsum{\langle\bfy\rangle \le N/\langle\bfw\rangle}{\bfy>\bfV} h(\bfw,\bfy)
$$
where $\bfV=(V^{(1)},\ldots, V^{(k-r)})$ runs through vectors with $V^{(j)}$ either $V_m$ or $V_{m+1}$ when $j\in{\cal L}_m$
with $m<k$,
and $V^{(j)} = V_k$ for $j\in{\cal L}_k$, and $\eps(\bfV)\in\{0,1\}$ is chosen 
appropriately. Note that $V^{(j)}\ge V_{l+1}$ for all $j$. Therefore, Lemma \ref{L2.9} and Lemma \ref{L2.3} 
yield
\be{2611}
\multsum{\langle\bfy\rangle \le N/\langle\bfw\rangle}{\bfy>\bfV} h(\bfw,\bfy)
= 
\f{c_{h,r}(\bfw)}{\langle\bfw\rangle^\al} N^\al p_{k-r} 
\Big(\al \log \f{N}{\langle\bfw\rangle \langle\bfV\rangle }\Big)
+ O(N^{\al+\eps} V_{l+1}^{-\delta/2k}  \langle\bfw\rangle^A) .
\ee
Now note that $\langle\bfV\rangle = N^\rho$ for some positive $\rho$, and hence that
$$ \log \f{N}{\langle\bfw\rangle \langle\bfV\rangle} = (1-\rho)\log N - \log 
\langle\bfw\rangle.
$$
Consequently, $p_{k-r} 
(\al \log {N}/{\langle\bfw\rangle \langle\bfV\rangle })$ can be written as a polynomial
in $ \log \langle\bfw\rangle $, with coefficients containing powers of $\log N$. It is then
possible to sum the equation \rf{2611} over $|\bfw|\le V_l$ by Theorem \ref{S2.7}, and the result 
over the finitely many $\bfV$. Since $\log V_l$ is a constant multiple of $\log N$,
one obtains a formula
$$  \Ups_{\fk A} (N) = N^\al P_{\fk A} (\log N) + O(N^{\al+\eps} V_{l+1}^{-\delta/2k} 
V_l^{r(A+1)}) $$
in which $P_{\fk A}$ is as desired, and an inspection of \rf{defV} shows that the error term
is indeed $O(N^{\al-\eta})$. This completes the proof of \rf{269} in all cases.

\def\la{\langle}
\def\ra{\rangle}

\section{Weyl sums over products}

\subsection{Introductory comment}
It is apparent that a circle method approach to count solutions of the equation \rf{01} will involve
the exponential sum
\be{321}
f(\al) = f_k(\al,\bfX) = \multsum{x_i\le X_i}{1\le i\le k} e(\al\la\bfx\ra^d)
\ee
in which $\bfX\in [1,\infty)^k$. The case $k=1$ is that of classical Weyl sums,
with an extensive literature. Little appears to be available for $k\ge 2$, 
forcing us to rework the most basic theory of Weyl sums in the new context.  
The simplest principles will be sufficient for our purposes. Most of our estimates may be improved,
but such refinements will not be needed here.  Throughout this section vectors are of dimension $k$.

\subsection{Approximate formulae}
For $q\in\NN$ and $\bfX$ as above let $E(q,\bfX)$ be the symmetric function
in $X_1,\ldots,X_k$ that, whenever
\be{322} X_1\ge X_2\ge \ldots \ge X_k\ee
holds, is defined by
\be{323} E(q,\bfX) = q^k + \sum_{r=1}^{k-1} q^{k-r}X_1\cdots X_r.\ee
We also define the complete Weyl sum
\be{324} S(q,a) = \multsum{x_i=1}{1\le i\le k}^q e(a\la\bfx\ra^d/q).\ee
Only in cases where  $k$ varies this  will be indicated 
by writing $E_k$ or $S_k$.

\begin{lem}\label{L321} Whenever $a\in\ZZ$, $q\in\NN$ and $\bfX\in[1,\infty)^k$, one has
$$ f_k(a/q,\bfX) = q^{-k} S(q,a) \la\bfX\ra + O(E(q,\bfX)).
$$
\end{lem}

\noindent
{\em Proof}. By symmetry, we may suppose that \rf{322} holds. Now sort the $x_j$
in \rf{321} into residue classes modulo $q$ to confirm that
$$ f\Big(\f{a}{q}\Big) = \sum_{\bfb=1}^q e\Big(\f{a\la\bfb\ra^d}{q}\Big)
            \multsum{\bfx\le \bfX}{\bfx\equiv\bfb \bmod q} 1 
= S(q,a)\prod_{j=1}^k \Big(\f{X_j}{q} + O(1)\Big).
$$
The product of all $X_j/q$ yields the leading term while an inspection of 
\rf{323} shows that all other terms are bounded by $E(q,\bfX)$.

\medskip

We now apply partial summation to evaluate $f(a/q+\beta)$. This features the function
\be{325}
v(\beta)=v_k(\beta,\bfX) = \int_0^{X_k} \cdots \int_0^{X_1} e(\beta\la\bft\ra^d)\,\dd t_1\ldots \dd t_k. \ee

\begin{lem}\label{L322} Whenever $a\in\ZZ$, $q\in\NN$, $\beta\in\RR$
 and $\bfX\in[1,\infty)^k$, one has
$$ f_k\Big(\f{a}q +\beta,\bfX\Big) = q^{-k} S(q,a)v(\beta)  + 
O\big(E(q,\bfX)(1+\la\bfX\ra^d |\beta|)^k \big).
$$
\end{lem}

{\em Proof}. Notation is the most difficult part of the otherwise  routine 
argument. Throughout this proof, let $\cal S$ denote a subset of $\{1,2,\ldots,k\}$, and write $\overline{\cal S} = \{1,\ldots,k\}\setminus\cal S$. As on earlier occassions, for $\bft\in[1,\infty)^k$, let $\bft_{\cal S} =(t_j)_{j\in\cal S}$. Also,
for $\bft,\bfu\in[1,\infty)^k$, let $(\bfu,\bft)_{\cal S}$ be the vector
$(z_1,\ldots,z_l)$ with
$$ z_j=t_j \quad \text{for }j\in{\cal S}, \quad 
z_j=u_j \quad \text{for } j\in\overline{\cal S}.$$
In the interest of brevity, we also write
$$ g(\bft) = e(\beta\la\bft\ra^d),$$
and if $\cal S$ consists of the $r$ numbers $s_1,\ldots,s_r$, then we put
$$ g^{(\cal S)}(\bft) = \f{\partial}{\partial t_{s_1}}\ldots
\f{\partial}{\partial t_{s_r}} g(\bft). $$

We are ready to apply partial summation to the sums over $x_j$ in \rf{321}.
This produces
\be{326}
f\Big(\f{a}q +\beta\Big) = g(\bfX)f\Big(\f{a}q\Big) + \sum_{{\cal S}\neq\emptyset}
(-1)^{\#\cal S} \int_{W({\cal S})} g^{(\cal S)}((\bfX,\bft)_{\cal S})
f\Big(\f{a}q, (\bfX,\bft)_{\cal S}\Big) \,\dd\bft_{\cal S}
\ee
in which the sum over $\cal S$ runs over subsets of $\{1,\ldots,k\}$,
and $W({\cal S})$ denotes the cartesian product of the invervals $[1,X_j)$
with $j\in\cal S$. The term $g(\bfX)f(a/q)$ may be considered as the formal
term ${\cal S }=\emptyset $ of the sum on the right.

We now apply Lemma \ref{L321} to all summands on the right hand side of
\rf{326}. The leading terms that arise reassemble to
\be{327}
q^{-k}S(q,a)\sum_{\cal S} (-1)^{\#\cal S} 
\int_{W({\cal S})} g^{(\cal S)}((\bfX,\bft)_{\cal S})
\la\bfX_{\overline{\cal S}}\ra \la\bft_{\cal S}\ra  \,d\bft_{\cal S},
\ee
with the summand ${\cal S}=\emptyset$ to be read as $g(\bfX)\la\bfX\ra$. 
If the formal partial integration 
$$ \int_0^X h(t)\,\dd t = Xh(X) - \int_0^X t h'(t)\,\dd t $$
is applied to all integrations in \rf{325}, then one finds that
$$ v(\beta) = g(\bfX)\la\bfX\ra + \sum_{{\cal S}\neq\emptyset}
(-1)^{\#\cal S}   \la\bfX_{\overline{\cal S}}\ra \int_{W({\cal S})} \la\bft_{\cal S}\ra 
g^{(\cal S)}((\bfX,\bft)_{\cal S})
\,\dd\bft_{\cal S}
.$$
Hence, the sum in \rf{327} is exactly the leading term on the right hand side of the formula in Lemma \ref{L322}. 

It remains to control the error terms that arise from the use of Lemma \ref{L321} in \rf{326}. The transition from $g(\bfX)f(a/q)$ to $q^{-k}S(q,a)\la\bfX\ra g(\bfX)$ results in an error bounded by $E(q,\bfX)$, which is acceptable. For the remaining terms, first note that for  $1\le r\le k$ there are natural numbers $b_{r,j}$
with
$$ \f{\partial}{\partial t_1} \ldots  \f{\partial}{\partial t_r} g(\bft)
= d^r g(\bft) \sum_{j=1}^r b_{r,j} (2\pi i\beta)^j \f{\la\bft\ra^{jd}}{t_1\cdots t_r}, 
$$
as one readily confirms by induction on $r$. Hence, by symmetry, it follows that whenever $\cal S$ is non-empty, one has
$$ g^{(\cal S)}(\bft) \ll \sum_{j=1}^{\#\cal S} |\beta|^j \f{\la\bft\ra^{jd}}{\la\bft_{\cal S}\ra}.$$
Moreover, we note that for $\bft_{\cal S}\in W({\cal S})$ one has
$$ E(q, (\bfX,\bft)_{\cal S}) \le E(q,\bfX).$$
Hence, the insertion of the asymptotic relation from Lemma \ref{L321} into \rf{326} is at the cost of an error not exceeding
\begin{eqnarray*}
& \ll & E(q,\bfX)\int_{W(\cal S)} |g^{(\cal S)}((\bfX,\bft)_{\cal S})|\, \dd\bft_{\cal S} \\
& \ll & E(q,\bfX)\sum_{j=1}^{\#\cal S} |\beta|^j \la\bfX_{\overline{\cal S}}\ra^{jd}\int_{W(\cal S)} \la\bft_{\cal S}\ra^{jd-1}\,\dd\bft_{\cal S} \\
& \ll & E(q,\bfX)\sum_{j=1}^{\#\cal S} |\beta|^j \la\bfX\ra^{jd}.
\end{eqnarray*}
Lemma \ref{L322} is now immediate.

\medskip

Before we continue our study of the exponential sum \rf{321}, we briefly estimate the factors $S(q,a)$ and $v(\beta)$ in the leading term of the asymptotic expansion supplied by Lemma \ref{L322}.

\begin{lem}\label{L323}
Let $q\in\NN$, $a\in\ZZ$ and $(a,q)=1$. Then
$$ q^{-k} S(q,a) \ll \tau(q)^{k-1} q^{-1/d}.$$
Moreover, when $a'\in\ZZ$, $q'\in\NN$ are any numbers with $a/q = a'/q'$, then
$$ q^{-k} S(q,a) = (q')^{-k} S(q',a'). $$
\end{lem}

\noindent
{\em Proof}. By two applications of Lemma \ref{L321}, we see that
$$  
q^{-k} S(q,a) = \lim_{X\to\infty} X^{-k} f_k(a/q, (X,\ldots,X)), $$
 and that this also holds  with $q',a'$ in place of $q,a$. However, the right hand side 
here remains the same if $a/q$ is replaced by $a'/q'$.
This already confirms the
second clause in Lemma \ref{L323}.

Turning our attention to the first clause, we note that the case $k=1$
is a familiar estimate of Hardy and Littlewood. More precisely, for $d\ge 2$,
it follows from Theorem 4.2 of \cite{hlm} and the preceding remark that
whenever $b\in\ZZ$ and $r\in\NN$, then
$$ \sum_{x=1}^r e(bx^d/r) \ll r^{1-1/d} (r,b)^{1/d}. $$
When $d=1$, orthogonality evaluates the sum on the left here, and 
 the bound remains valid. Hence, for $k=1$, the proof of the lemma is complete,
and for $k\ge 2$ and $(a,q)=1$, we have
\begin{eqnarray*}
q^{-k}S_k(q,a) &=& \sum_{x_2,\ldots,x_k=1}^q q^{-k} S_1(q,a(x_2\ldots x_k)^d)\\
&\ll& q^{1-k-1/d} \sum_{x_2,\ldots,x_l=1}^q (q,(x_2\ldots x_k)^d)^{1/d}\\
& \ll& q^{-1/d} \Big(\sum_{x=1}^q \frac{(q,x)}{q}\Big)^{k-1}.
\end{eqnarray*}
The desired estimate is now immediate. 

\medskip

We close this section with an upper bound for the integral $v_k(\beta,\bfX)$
introduced in \rf{325}. Let
\be{329}
V_k(\beta) = \int_{[0,1]^k} e(\beta\la\bft\ra^d)\,\dd\bft
\ee
which is the special case $X_j=1$ $(1\le j\le k)$ of \rf{325}. An obvious substitution yields
\be{3210}
v_k(\beta,\bfX) = \la\bfX\ra V_k(\la\bfX\ra^d\beta). 
\ee
By \rf{329}, one finds that
\be{3211} V_k(\beta)\ll 1. 
\ee
We now proceed to show that whenever $|\beta|\ge 1$, then
\be{3212}
V_k(\beta) \ll |\beta|^{-1/d} (1+ \log |\beta|)^{k-1}.
\ee
In fact, the bound $V_1(\beta) \ll |\beta|^{-1/d}$ is immediate by partial
integration (or explicit computation of the integral when $d=1$). We now proceed by induction on $k$. 
By \rf{329},
\be{3213}
V_{k+1}(\beta) = \int_0^1 V_k(\beta t^d) \,\dd t.
\ee
In the range $0\le t\le |\beta|^{-1/d}$ we use \rf{3211} to see that these $t$
contribute at most $O(|\beta|^{-1/d})$ to the integral in \rf{3213}. With
\rf{3212} at our disposal, it follows that
$$ V_{k+1}(\beta) \ll |\beta|^{-1/d} + \int_{|\beta|^{-1/d}}^1 |\beta|^{-1/d}t^{-1} (1+\log (|\beta|t^d))^{k-1}\,\dd t. $$
This implies \rf{3212} with $k+1$ in place of $k$, completing the induction.

We may summarise the conclusions in \rf{3210}, \rf{3211} and \rf{3212}
in the following lemma.

\begin{lem}\label{L324} Let $\bfX\in[1,\infty)^k$ and $\beta\in\RR$. Then
$  v_k(\beta,\bfX) \ll \la\bfX\ra ,$
and whenever $|\beta| \ge \la\bfX\ra^{-d}$, one has
$$  v_k(\beta,\bfX) \ll |\beta|^{-1/d} \big(1+ \log (\la\bfX\ra^{d}|\beta|)\big)^{k-1}.$$
\end{lem}

\subsection{Bounds of Weyl's type}
Our next goal is a crude form of Weyl's inequality for the sum \rf{321}. This will then be coupled with the results of the previous section
to provide suitable pointwise minor arc estimates. Throughout this section, we write
$$ D= 2^{d-1}. $$

\begin{lem}\label{L331}
Suppose that $\bfX\in[1,\infty)^k$ satisfies \rf{322}. Let $\al\in\RR$, $a\in\ZZ$, $q\in\NN$
with $(a,q)=1$ and $|q\al -a|\le 1/q$. Then
$$ |f(\al)|^{D^k} \ll \la\bfX\ra^{D^k+\eps} \Big(\f1q + \f1{X_k} + \f{q}{\la\bfX\ra^d} \Big). $$
\end{lem}

{\em Proof}. First suppose that $d=1$. Then $D=1$. We carry out the sum over $x_1$. Then, by a familiar divisor function estimate,
\begin{eqnarray*}
f(\al) &\ll & \multsum{x_j\le X_j}{2\le j\le k} \min(X_1,\parallel\al x_2\cdots x_l\parallel^{-1}) \\
& \ll & \la\bfX\ra^\eps \sum_{u\le X_2\cdots X_k} \min (X_1,\parallel \al u\parallel^{-1}).
\end{eqnarray*}
The desired estimate now follows from Lemma 2.2 of Vaughan \cite{hlm}.

\smallskip

Next suppose that $d\ge 2$ and recall that Weyl's differencing argument produces the inequality
\be{331}
\Big| \sum_{y\le Y} e(\gamma y^d)\Big|^D \le (2Y)^{D-d} \sum_{|\bfh|< Y}
\sum_{y\in I(\bfh)} e(d! \la\bfh\ra y\gamma)
\ee
in which $\bfh=(h_1,\ldots,h_{d-1})$ and $I(\bfh)$ is the set of all
$y\in\NN$ with
$$ 1\le y\le Y, \quad 1\le y+h_j \le Y \quad (1\le j\le d-1)$$
(see Vaughan \cite{hlm}, Lemma 2.3 and Exercise 2.8.1).
We now prove by induction on $k$ that
\be{332}
|f_k(\al,\bfX)|^{D^k} \le (2^k \la\bfX\ra)^{D^k-d} 
\multsum{|\bfh_1|< X_1}{x_1\in I_1(\bfh_1)} \cdots
\multsum{|\bfh_l|< X_l}{x_k\in I_k(\bfh_k)} 
e\Big( (d!)^k \al \prod_{j=1}^k \la\bfh_j\ra x_j\Big)
\ee
where $\bfh_j \in\ZZ^{d-1}$, and $I_j(\bfh_j)$ is a certain interval contained
in $[1,X_j]$. In fact, the case $k=1$ is \rf{331} with $Y=X_1$, $\gamma=\al$.
If $k\ge 2$ and \rf{332} is already confirmed for $k-1$ in place of $k$,
then one first applies H\" older's inequality to the identity
$$ f_k(\al,\bfX) = \sum_{x_k\le X_k} f_{k-1}(\al x_k^d, (X_1,\ldots,X_{k-1}))
$$
to infer that
$$ |f_k(\al,\bfX)|^{D^k} \le X_k^{D^k-D} 
\Big( \sum_{x_k\le X_k} |f_{k-1} (\al x_k^d, (X_1,\ldots,X_{k-1}))|^{D^{k-1}}\Big)^D, 
$$
and is then in a position to apply the induction hypothesis to
$ f_{k-1} (\al x_k^d, (X_1,\ldots,X_{k-1}))$. This yields
$$
|f_k(\al,\bfX)|^{D^k} \le X_k^{D^k-D} (2^{k-1}X_1\cdots X_{k-1})^{D^k-Dd}
\Big( \trisum{|\bfh_j|< X_j}{x_j\in I_j(\bfh_j)}{1\le j\le k-1} \sum_{x_k\le X_k} e(\al\Phi)\Big)^D $$
with
$$ \Phi = (d!)^{k-1} x_k^d \prod_{j=1}^{k-1} \la\bfh_j\ra x_j.$$
Here, the innermost sum over $x_k$ is again an ordinary $d$-th power Weyl sum.
By H\"older's inequality again, we first deduce that
$$
|f_k(\al,\bfX)|^{D^k} \le X_k^{D^k-D} (2^{k-1}X_1\cdots X_{k-1})^{D^k-d}
 \trisum{|\bfh_j|< X_j}{x_j\in I_j(\bfh_j)}{1\le j\le k-1} 
\Big| \sum_{x_k\le X_k} e(\al\Phi)\Big|^D, $$
and may then apply \rf{331} to arrive at \rf{332}, thus completing the induction.

Now consider the product
$$u = (d!)^{k} \prod_{j=1}^k \la\bfh_j\ra x_j $$
that occurs in \rf{332}. With $\bfh_j$, $x_j$ subject to the conditions of summation in \rf{332}, one has $u=0$ if and only if one of the components of some
$\bfh_j$ vanishes. Hence, by \rf{322}, the number of such $\bfh_j$, $x_j$
where $u=0$ is bounded by $\la\bfX\ra^d X_k^{-1}$, and the contribution of these
terms to \rf{332} amounts to $O(\la\bfX\ra^{D^k} X_k^{-1})$ which is acceptable.
For the remaining terms, we write $u=x_1m$ and sum over $x_1$ first. A divisor
function estimate then delivers the bound
 $$
|f_k(\al,\bfX)|^{D^k} \ll \la\bfX\ra^{D^k} X_k^{-1} + \la\bfX\ra^{D^k+\eps-d}
\sum_{1\le m\le M} \min (X_1, \parallel \al m\parallel^{-1})
$$
in which $M=(d!)^k \la\bfX\ra^d X_1^{-1}$. The conclusion of Lemma \ref{L331}
is now immediate from Lemma 2.2 of \cite{hlm}.

\medskip

In preparation for two applications of the Hardy-Littlewood method, we introduce a dissection into major and minor arcs. From now on, we write
\be{defP} P = \la\bfX\ra. \ee
With $Q$ in the range $1\le Q\le P^{1/(dk)}$, let ${\fk N}(Q)$ denote the
union of the pairwise disjoint intervals
$|q\al-a|\le QP^{-d}$ with $a,q$ subject to $1\le q\le Q$, $a\in\ZZ$ 
and $(a,q)=1$. Let ${\fk n}(Q)= \RR\setminus {\fk N}(Q)$ and define
${\fk M}(Q)$ and ${\fk m}(Q)$ by
$$ {\fk M}(Q) = {\fk N}(Q)\cap [0,1],\quad {\fk m}(Q) = {\fk n}(Q)\cap [0,1].$$

\begin{lem} \label{weyl1}  Given natural numbers $d,k$, let 
$$ \om= (8dk)^{-8}, \quad \eta=\om/(2D^k dk).$$
Then
$$ \sup_{\al \in{\fk n}(P^\om)} |f_k(\al,\bfX)| \ll P^{1-\eta}. $$
\end{lem}

\noindent
{\em Proof}. By symmetry, we may suppose that \rf{322} holds.
There are two cases, depending on the relative size
of the $X_j$. 

First suppose that $X_k \ge P^{\om/(dk)}$. By Dirichlet's theorem, choose
coprime integers $a,q$ with $1\le q\le P^{d-\om}$ and $|q\al -a|\le P^{\om-d}$.
Since $\al\in{\fk n}(P^{\om})$, it follows that $q>P^\om$. Lemma \ref{L331}
now yields $f_k(\al,\bfX)\ll P^{1-2\eta+\eps}$.  

It remains to consider the case where $X_k< P^{\om/(dk)}$. By \rf{322},
one has $X_1\ge P^{1/k}$, and hence, there is a number $r$ with $1\le r\le k-1$
and
\be{345}
X_r\ge P^{\om/(dk)}> X_{r+1}.
\ee
We write $\bfY=(X_1,\ldots,X_r)$ and $\bfZ=(X_{r+1},\ldots,X_k)$. Then, by
\rf{321},
\be{346}
f_k(\al,\bfX) = \sum_{\bfz\le \bfZ} f_r(\al\la\bfz\ra^d, \bfY).
\ee

Now consider $\bfz$ with $\bfz\le \bfZ$. By Dirichlet's theorem,
there are coprime $a=a(\bfz)$, $q=q(\bfz)$ with $1\le q\le P^{d-\om}$
and $|q\la\bfz\ra^k\al -a |\le P^{\om-k}$. Note that \rf{345}
implies 
\be{Zet} \la\bfz\ra\le \la\bfZ\ra \le P^{\om(k-1)/(dk)}. 
\ee
Hence, for $q\le P^{\om/k}$, it follows that $q\la\bfz\ra^d \le P^\om$,
so that $\al\in{\fk N}(P^\om)$. This is not the case. Consequently, we may conclude that we must have $q>P^{\om/k}$, and Lemma \ref{L331} yields
\begin{eqnarray*}
|f_r(\al\la\bfz\ra^d,\bfY)|^{D^r} &\ll &
\la\bfY\ra^{D^r+\eps} \Big( \f1q + \f1{X_r} + \f{q}{\la\bfY\ra^d}\Big)\\
&\ll&   \la\bfY\ra^{D^r+\eps} (P^{-\om/(dk)} + \la \bfZ\ra^d P^{-\om})\\
&\ll& \la\bfY\ra^{D^r+\eps} P^{-\om/(dk)}.
\end{eqnarray*}
In the last line, we have used \rf{Zet}. We take the $D^r$-th root and then
sum over $\bfz$. By \rf{346}, it follows that $f_k(\al,\bfX)\ll P^{1-\eta}$. 
This completes the proof.

\begin{lem}\label{weyl2} 
Let $d,k,\om,\eta$ be as in the previous lemma, and let $s$ be a natural number.
Suppose that \rf{322} holds. Put
$ U=  X_k^{k\om/s}. $
Then
$$ \sup_{\al\in{\fk n}(U)}|f_k(\al,\bfX)|\ll P^{1-\eta} + PU^{-1/(5dk)}.$$
\end{lem}

\noindent
{\em Proof}. 
Note that $X_k\le P^{1/k}$ so that $U\le P^\om$.  Whenever $\al\in {\fk n}(P^\om)$, then
Lemma \ref{weyl1}
 supplies a satisfactory bound. Hence, we may concentrate on the case
where $\al\in {\fk N}(P^\om)\cap {\fk n}(U)$.

Recall that $X_1\ge P^{1/k}$. Also, note that $3k^4\om < 1/k$. Consequently,
in the current situation, there is a number $r$ with $1\le r\le k-1$ and
$$ X_r\ge P^{3k^4\om} >X_{r+1}.$$ 
Note that this implies $k\ge 2$.
With this new definition of $r$, we use the identity \rf{346} from the previous
proof. The estimation now proceeds through the inequality
\be{fdec} |f(\al)| \le \sum_{\bfz\le \bfZ} |f_r(\al\la\bfz\ra^d, \bfY)|, \ee
where we use the results from the preceding section to bound
$f_r(\al\la\bfz\ra^d,\bfY)$. Since $\al\in{\fk N}(P^\om)$, there is a unique
pair $a,q$ with $1\le q\le P^\om$, $a\in\ZZ$, $(a,q)=1$ and $\al = (a/q) + \beta$ satisfying $|q\beta|\le P^{\om-d}$. We use Lemma \ref{L322} with
$r$ in place of $k$, and with $a$ and $\al$ replaced by $a\la\bfz\ra^d$,
 $\al\la\bfz\ra^d$. Then, by a crude use of Lemmas \ref{L323} and \ref{L324},
one finds that
$$ f_r(\al\la\bfz\ra^d,\bfY) \ll (q,\la\bfz\ra^d)^{1/d} \la\bfY\ra
(q+ q\la\bfY\ra^d \la\bfz\ra^d |\beta|)^{\eps-1/d} + E(q,\bfY)
(1+ \la\bfY\ra^d \la\bfz\ra^d |\beta|)^r.$$
Here, we first concentrate on the second term on the right hand side.
For $\bfz\le\bfZ$, the upper bound on $|\beta|$ guarantees that
$$ 1+ \la\bfY\ra^d \la\bfz\ra^d |\beta| \le 1 + P^d|\beta| \le 2P^\om q^{-1}.$$
Moreover, in the current context, we have $q\le P^\om\le X_j$ for all $j\le r$,
so that \rf{323} now shows
$$ E(q,\bfY) (1+ \la\bfY\ra^d \la\bfz\ra^d |\beta|)^r
\le 2^r P^{\om r} \sum_{j=0}^{r-1} q^{-j} \prod_{i=1}^j X_i .$$
Here, the right hand side does not exceed $O(\la\bfY\ra X_r^{-1} P^{\om r})$,
and by the definition of $r$ this is bounded by $O(\la\bfY\ra P^{-\om})$. 
We now insert the results obtained so far into \rf{fdec} and find that
\begin{eqnarray*} 
f(\al) &\ll& P^{1-\om} + \la\bfY\ra q^{\eps-1/d} \sum_{\bfz\le\bfZ}
(q,\la\bfz\ra) (1+ \la\bfY\ra^d\la\bfz\ra^d |\beta|)^{\eps-1/d}\\
&\ll& P^{1-\om} + \min \Big( \la\bfY\ra q^{\eps-1/d}\sum_{\bfz\le\bfZ}
(q,\la\bfz\ra),\,  \la\bfY\ra^{d\eps} |\beta|^{\eps-1/d}\sum_{\bfz\le\bfZ}
(q,\la\bfz\ra) \la\bfz\ra^{d\eps-1} \Big).
\end{eqnarray*}
\smallskip
The simple estimate
\be{qeps} \sum_{\bfz\le\bfZ} (q,\la\bfz\ra)
\le \prod_{j=r+1}^k \sum_{x_j\le X_j} (q,x_j) \ll q^\eps \la\bfZ\ra \ee
and partial summation now show that
$$ f(\al ) \ll P^{1-\om} + P q^{\eps-1/d} + P^{d\eps} |q\beta|^{\eps-1/d} $$
For  $\al\in{\fk n}(U)$, we have $q>U$ or $q|\beta|> UP^{-d}$, and the
lemma follows.

\subsection{A mean value estimate}
The mean value \rf{03} defining $n_0(d)$ may play the role of Hua's lemma in
an investigation of diagonal forms via the circle method. We need a similar bound
for Weyl sums over products. In fact, because  the number of solutions
of $\la\bfx\ra =m$ with $\bfx\in\NN^k$ is $O(m^\eps)$, a consideration of
the underlying diophantine equations shows that
\be{341}
\int_0^1 |f_k(\al,\bfX)|^{n_0(d)} \, d\al\ll \la\bfX\ra^\eps \int_0^1 |f_1(\al,
\la\bfX\ra)|^{n_0(d)}\,d\al   \ll \la\bfX\ra^{n_0(d)-d+2\eps}.
\ee
In cases where some of the $X_j$ are considerably smaller than $\la\bfX\ra^\eps$
this estimate is insufficient for our purposes. The following lemma provides an
alternative bound, free of the unwanted $\eps$ in \rf{341}.

\begin{lem}\label{L341}
Fix a real number $\sgm> n_0(d)$. Then, for $\bfX\in [1,\infty)^k$ one has
$$ \int_0^1 |f_k(\al,\bfX)|^\sgm\,d\al \ll \la\bfX\ra^{\sgm-d}. $$
\end{lem}  

{\em Proof}. By symmetry, we may suppose that \rf{322} holds. Recall that
this implies $X_1\ge P^{1/k}$ where $P$ is defined by \rf{defP}. Also, let
$\om$ be as in Lemma \ref{weyl1}. Now, by that lemma and \rf{341},
one has
\be{mino}
\int_{{\fk m}(P^\om)} |f(\al)|^\sgm\,d\al
\le \sup_{\al\in{\fk m}(P^\om)} |f(\al)|^{\sgm-n_0} \int_0^1 |f(\al)|^{n_0}\,d\al \ll P^{\sgm-d}. 
\ee

It remains to consider the major arcs ${\fk M}(P^\om)$ where one uses \rf{346}
with $r=1$ and $\bfZ=(X_2,\ldots,X_k)$. If $\al\in{\fk M}(P^\om)$, let $a,q$ be
the unique pair $a,q$ with $1\le q\le P^\om$, $(a,q)=1$ and $\al=a/q+\beta$
satisfying $q|\beta|\le P^{\om-d}$. Then, by Lemma \ref{L322},
$$ f(\al) = \sum_{\bfz\le \bfZ } \big(
q^{-1}S_1(q,a\la\bfz\ra^d)v_1(\beta\la\bfz\ra^d,X_1) +O\big(q(1+X_1^d\la\bfz\ra^d|\beta|)\big)\big),$$
and for $\al\in{\fk M}(P^\om)$ one has
$q(1+X_1^d\la\bfz\ra^d|\beta|) \ll P^\om$.
Hence, on writing
$$ H(q,\beta) = \sum_{\bfz\le \bfZ }(q,\la\bfz\ra) \min(X_1, \la\bfz\ra^{-1}|\beta|^{-1/d}), $$
one now infers from Lemmas \ref {L323} and \ref{L324} that
\be{B1} f(\al) \ll q^{-1/d} H(q,\beta) + \la\bfZ\ra P^\om  \ll   q^{-1/d}H(q, \beta) + P^{1 - 1/(2k)}, \ee
as is apparent from the  lower bound on $X_1$. 
By \rf{qeps}, 
$$ H(q,\beta)\le X_1 \sum_{\bfz\le \bfZ} (q,\la\bfz\ra) \ll Pq^\eps. $$
The previous bound for $f$ now delivers
\be{B2} f(\al) \ll Pq^{\eps-1/d} + P^{1-1/(2k)} \ll  Pq^{-1/(2d)} \ee
uniformly for $\al\in{\fk M}(P^\om)$.
Since $n_0(d)\ge 2d$, we may use \rf{B1} for $2d$ copies of $f$ and \rf{B2} for the rest  to confirm that
$$ |f(\al)|^\sgm \ll (Pq^{-1/2d})^{\sgm-2d} (q^{-2}H(q,\beta)^{2d} + P^{2d-(d/k)}). $$
The measure of ${\fk M}(P^\om)$ is $O(P^{2\om-d})$. Hence, on integrating the previous inequality, we infer that
$$ \int_{{\fk M}(P^\om)} |f(\al)|^\sgm\,d\al \ll P^{\sgm-2d}I + P^{\sgm-d-1/(2k)} $$
where
$$ I =  \sum_{q\le P^\om} q^{-\sgm/(2d)} \int_{-\infty}^\infty H(q,\beta)^{2d}\,d\beta. $$
Now
$$ \int_{-\infty}^\infty H(q,\beta)^{2d}\,d\beta
= \multsum{\bfz_j\le \bfZ}{1\le j\le 2d} 
(q,\la\bfz_1\ra)\cdots (q,\la\bfz_{2d}\ra) 
\int_{-\infty}^\infty \prod_{j=1}^{2d} \min(X_1, \la\bfz\ra^{-1}|\beta|^{-1/d})\,d\beta. $$
By H\"older's inequality and \rf{qeps}, it follows that
\begin{eqnarray*}
\int_{-\infty}^\infty H(q,\beta)^{2d}\,d\beta 
& \le & \Big( \sum_{\bfz\le\bfZ} (q,\la\bfz\ra) \Big(\int_{-\infty}^\infty
\min(X_1^{2d}, \la\bfz\ra^{-2d} |\beta|^{-2})\,d\beta \Big)^{1/(2d)} \Big)^{2d} \\
& \ll & \Big( \sum_{\bfz\le\bfZ} (q,\la\bfz\ra) X_1^{1/2}\la\bfz\ra^{-1/2} \Big)^{2d} 
\\ &\ll& (q^\eps X_1^{1/2}\la\bfZ\ra^{1/2})^{2d} \ll P^dq^{2d\eps}.
\end{eqnarray*}
Since $\sgm>2d$, one now first confirms that $I\ll P^d$, and then that
$$  \int_{{\fk M}(P^\om)} |f(\al)|^\sgm\,d\al \ll P^{\sgm-d}. $$
The lemma follows on combing this bound with
\rf{mino}.

\section{Multihomogeneous diagonal forms}

\subsection{The auxiliary theorem}

In the present chapter we consider integer solutions of equations
similar to \rf{01}. 
Fix natural numbers $d,k,s$. 
Let $c_j\in\ZZ\setminus\{0\}$ for $1\le j\le s$. With $\bfx_j\in\ZZ^k $, consider
the diophantine equation
\be{31}
\sum_{j=1}^s c_j\la\bfx_j\ra^d =0.\ee
For $\bfX\in[1,\infty)^k$, let $\MM_\bfc(\bfX)$ denote the number of
 solutions of \rf{31} with
$$ 1\le |x_{j,r}| \le X_r \quad (1\le j\le s,\,1\le r\le k), $$
and let $\MM^+_\bfc(\bfX)$ be the number of the solutions with all $x_{j,r}$ positive.

\begin{thm}\label{S4.1} Let $s>n_0(d)$. Then, there is a 
positive number $\delta$ such that whenever $c_j\in\ZZ\setminus\{0\}$, 
there  are non-negative real numbers ${\fk E}(\bfc)$,  ${\fk E}^+(\bfc)$ 
with
\be{32} \MM_\bfc(\bfX) = {\fk E}(\bfc)\la\bfX\ra^{s-d}
+O(\la\bfX\ra^{s-d} (\min X_i)^{-\delta} |\bfc|^{s+k}) \ee
and
\be{32+} \MM_\bfc^+(\bfX) = {\fk E}^+(\bfc)\la\bfX\ra^{s-d}
+O(\la\bfX\ra^{s-k} (\min X_i)^{-\delta} |\bfc|^{s+k}).\ee
Further, the number ${\fk E}(\bfc)$ is positive 
if and only if the equation
\be{33} c_1y_1^d + \ldots + c_s y_s^d =0 \ee
admits  non-trivial solutions in $\RR$ and in $\QQ_p$ for all primes $p$.
Also,  ${\fk E}^+(\bfc)$ is positive if and only if \rf{33} has non-trivial
solutions in  $\QQ_p$ for all primes $p$, and in positive real numbers.
\end{thm}

Note that we claim the asymptotic formulae \rf{32} and \rf{32+} with error terms that are uniform in $\bfc$.
However, we will make no attempt to estimate the leading term from below in terms of $\bfc$, as this
is not required later.

We shall prove \rf{32+} by a straightforward application of the circle method
in Section 4.3 below, and then deduce \rf{32} by a combinatorial observation in Section 4.4.
As is to be expected,
the constants  ${\fk E}(\bfc)$,  ${\fk E}^+(\bfc)$ turn out to be the product of the singular series
and the singular integral associated with the respective counting problem. These natural interpretations of
 ${\fk E}(\bfc)$ and  ${\fk E}^+(\bfc)$ can be used to establish our claims concerning positivity of these numbers.
We prepare a swift treatment of the local part of our circle method work
with a discussion of this matter in section 4.2. Once Theorem \ref{S4.1} is established, we close this chapter with a brief treatment of primitive solutions counted by $\MM_\bfc(\bfX)$.

\subsection{Singular series and integral} 
In preparation for the circle method work in the following section, we define and estimate the singular series
and the singular integral for the equation \rf{31}. We use the notation introduced in section 3.2. 
Fix a natural number $k$.
Throughout this section, let $c_1,\ldots,c_s$
be non-zero integers. Now recall \rf{324} and put
$$ T_\bfc(q) = q^{-ks} \multsum{a=1}{(a;q)=1}^q \prod_{j=1}^s S_k(q,ac_j) . $$
Then, by 
Lemma \ref{L323},
$$  T_\bfc(q) \ll q^{1-(s/d)+\eps} (q,c_1)^{1/d} \cdots (q,c_s)^{1/d} \ll q^{1-(s/d)+\eps} |\bfc|^{s/d}. $$
It follows that whenever $s>2d$, the singular series 
\be{singser}{\fk S}(\bfc)= {\fk S}_k(\bfc) = \sum_{q=1}^\infty T_\bfc(q)\ee
converges absolutely, and one has
\be{44} {\fk S}(\bfc) \ll |\bfc|^{s/d}.\ee
Moreover, for the partial sum
\be{45} {\fk S}(\bfc, W)=\sum_{q\le W} T_\bfc(q) \ee
we also have
\be{46} {\fk S}(\bfc, W)= {\fk S}(\bfc) + O(W^{-1/(2d)} |\bfc|^{s/d}). \ee

These immediate estimates already suffice within the analysis to be performed in the next section.
We now show that $ {\fk S}(\bfc)$ is a product of local densities. In fact,
by the argument underpinning the proof of Lemma 2.11 in Vaughan \cite{hlm}, one finds that
$T_\bfc(q)$ is multiplicative in $q$. Hence,  ${\fk S}(\bfc)$ can be written as an Euler product,
and by a suitable analogue of Lemma 2.12 of Vaughan \cite{hlm}, one may compute the Euler factors.
This yields
\be{47} {\fk S}(\bfc) = \prod_p E_p(\bfc) \ee
where
\be{Ep} E_p(\bfc) = \sum_{l=0}^\infty T_\bfc (p^l) = \lim_{L\to\infty} p^{L(1-ks)} \Phi_\bfc(p^L),\ee
and in which $\Phi_\bfc(q)$ denotes the number of incongruent solutions to
\be{Cong} \sum_{j=1}^s c_j\la\bfx_j\ra^d \equiv 0\bmod q. \ee
In particular, the number ${\fk S}(\bfc)$ is non-negative.

Now suppose that the equation \rf{33} has a non-trivial solution in $\QQ_p$. Since \rf{33} is homogeneous,
there will be a solution $\bfy\in\ZZ_p^s$ with $p\nmid \bfy$. By symmetry, we may suppose that $p\nmid y_1$.
The theory of $d$-th power residues supplies a natural number $\gamma$ with the property that for all integers $b$
for which the congruence $c_1z^d\equiv b \bmod p^\gamma$ has a solution with $p\nmid z$, the allied congruences  
 $c_1z^d\equiv b \bmod p^L$ are also soluble with $p\nmid z$, for all $L>\gamma$. Now choose integers $z_j$
with $z_j\equiv y_j \bmod p^\gamma$. 
Then $p\nmid z_1$. Now put $\bfz_j=(z_j,1,1,\ldots,1)$. Then, observing that $z_j=\la\bfz_j\ra$, we infer
$$ c_1\la\bfz_1\ra^d + \ldots+ c_s \la\bfz_s\ra^d \equiv 0 \bmod p^\gamma. $$
Let $L>\gamma$, and choose $\bfx_2,\ldots,\bfx_s$ modulo $p^L$ with $\bfx_j\equiv \bfz_j \bmod p^\gamma$.
There are $p^{(L-\gamma)k(s-1)}$ choices for $\bfx_2,\ldots,\bfx_s$ that are incongruent modulo $p^L$. By construction,
$$ c_1z_1^d \equiv -c_2\la\bfx_2\ra^d - \ldots - c_s\la\bfx_s\ra^d \bmod p^\gamma, $$
and hence there is a number $z$ with $p\nmid z$ and 
$$ c_1z^d \equiv -c_2\la\bfx_2\ra^d - \ldots - c_s\la\bfx_s\ra^d \bmod p^L, $$
It follows that any $\bfx_1$ with $\la\bfx_1\ra \equiv z \bmod p^L$ yields a solution of \rf{Cong},
and there are $\varphi(p^L)^{k-1}$ such $\bfx_1$ that are incongruent, modulo $p^L$. This proves that
$$ \Phi_\bfc (p^L) \ge \varphi(p^L)^{k-1} p^{(L-\gamma)k(s-1)}. $$
By \rf{Ep}, we conclude that $E_p(\bfc)>0$ whenever \rf{33} has non-trivial solutions in $\QQ_p$.
Since the product in \rf{47} converges absolutely, we may conclude as follows.

\begin{lem}\label{Lss} Let $s> 2d$. Then, for all natural numbers $k$, the singular series \rf{singser} is 
real and non-negative. If the equation \rf{33} admits non-trivial solutions in $\QQ_p$ for all primes $p$,
then the singular series is positive.
\end{lem}

We now turn to the singular integral. Note that the function $ t^{-1/d} (1+\log t)$ 
is decreasing on the interval  $t\ge 3$ . Hence when $|\beta|\ge 3$, we infer from \rf{3212} that the function
\be{48} {\fk V}_\bfc (\beta) = V_k(c_1\beta) V_k(c_2\beta) \ldots V_k(c_s\beta) \ee
obeys the inequality
$$  {\fk V}_\bfc (\beta) \ll |\beta|^{-s/d} (1+ \log |\beta|)^{s(k-1)} $$
uniformly in $\bfc$. Similarly, \rf{3211} gives $ {\fk V}_\bfc (\beta) \ll 1$. Hence, for $s>2d$, the 
singular integral
\be{49} {\fk I}^+(\bfc) = {\fk I}^+_k(\bfc) = \int_{-\infty}^\infty  {\fk V}_\bfc (\beta)\, \dd\beta \ee
exists, and one has
\be{410}  {\fk I}^+(\bfc) \ll 1 \ee
uniformly with respect to $\bfc$. By \rf{48}, we also see that whenever $s>2d$ and $W\ge 3$ one has
$$ \int_{W}^\infty  |{\fk V}_\bfc (\beta)|\, \dd\beta\ll W^{-1}. $$
Thus, the truncated singular integral
\be{411} {\fk I}^+(\bfc,W) =  \int_{-W}^W  {\fk V}_\bfc (\beta)\, \dd\beta
\ee
compares to \rf{49} through
\be{412}
{\fk I}^+(\bfc,W) = {\fk I}^+(\bfc) + O(W^{-1}).
\ee

Again, these simple estimates will suffice for our purposes. The arithmetical nature of $ {\fk I}^+(\bfc)$ is the subject of the next lemma.

\begin{lem}\label{Lsi}  Let $s> 2d$. Then, for all natural numbers $k$, the singular integral \rf{49} is 
real and non-negative. If the  $c_j$ are not all of the same sign,
then the singular integral is positive.
\end{lem}

{\em Proof}.
 The classical case $k=1$ has been worked out by Davenport \cite{D}, chapters 4 and 10. Although the set-up there is slightly different from ours, Davenport's argument immediately delivers the conclusions announced in Lemma \ref{Lsi}.
Hence, we now suppose that $k\ge 2$ and take up the story at \rf{329} where we substitute $\nu = t_1^dt_2^d\cdots t_s^d$
for $t_1$. On writing $\bft' =(t_2,\ldots,t_s)$, we then have
\be{413} V_k(\beta) = \f{1}{d} \int_{[0,1]^{k-1}} \la\bft'\ra^{-1} \int_0^{\la\bft'\ra^d} \nu^{(1/d)-1} e(\beta\nu)\,\dd\nu\,\dd\bft'. \ee
We may exchange the order of integration, as is most readily justified by Tonelli's theorem. For a compact presentation of the
outcome of this manouvre, let $\nu>0$ and put
$$ {\cal U}(\nu) = \{\bft'\in[0,1]^{k-1}: \,\la\bft'\ra ^d \ge \nu\}, \quad \psi_k(\nu) =  \nu^{(1/d)-1} \int_{{\cal U}(\nu)} \f{\dd\bft'}{\la\bft'\ra} . $$
Note that $\psi_k$ is integrable over $(0,1]$, as one readily confirms. By \rf{413}, we now deduce that
$$ V_k(\beta) = \f{1}{d} \int_0^1 \psi_k(\nu)e(\beta\nu)\,\dd\nu. $$
By \rf{48} and \rf{49},
$$ {\fk I}^+(\bfc) = d^{-s} \int_{-\infty}^\infty \int_{[0,1]^s} \psi_k(\nu_1)\ldots\psi_k(\nu_s)
e(\beta\bfc\cdot\bm\nu)\,\dd\bm\nu\,\dd\beta $$
where $\bfc\cdot\bm\nu= c_1\nu_1+\ldots+c_s\nu_s$ is the standard inner product. Define $t$ through the 
equation $c_1t = \bfc\cdot\bm\nu$ and substitute $t$ for $\nu_1$ in the inner integral. This produces the identity
\be{419}  {\fk I}^+(\bfc) = d^{-s} \int_{-\infty}^\infty \int_{-\infty}^\infty B(t) e(c_1  \beta t)\,\dd t\,\dd\beta \ee
where 
\be{420} B(t) =\int_{{\cal B}(t)} \psi_k(\nu_2)\ldots\psi_k(\nu_s)\psi_k\Big(t- \frac{\bfc'\cdot\bm\nu'}{c_1}\Big)\,\dd\bm\nu', \ee
and therein we put
\be{421} {\cal B}(t) =\Big\{ \bm\nu'\in [0,1]^{s-1}:\, 0\le t -  \frac{\bfc'\cdot\bm\nu'}{c_1}\le 1 \Big\}. \ee 
Note that $B(t)$ is a compactly supported continuous function, whence by \rf{419} and Fourier's integral theorem,
we infer that ${\fk I}^+(\bfc)= |c_1|^{-1}d^{-s} B(0)$. Since the integrand in \rf{420} is non-negative, this implies that
 ${\fk I}^+(\bfc)\ge 0$. Further, note that $s>2d$ implies that $s\ge 3$. Hence, if $c_1,\ldots,c_s$ are not all of the same sign, we can permute the indices to arrange that $c_2$ and $c_3$ have opposite signs, say. It is then immediate from \rf{421} that
${\cal B}(0)$ contains a box of positive $s-1$-dimensional volume on which the integrand in \rf{420}
is continuous and positive, whence $B(0)>0$. This completes the proof of the lemma.

\medskip

We close this section with a remark concerning the product $ {\fk I}^+(\bfc) {\fk S}(\bfc)$.
Suppose that that $\bfc$ is of the form $\bfc= a \bfc'$ with $a\in\NN$ and $\bfc'$ primitive.
Then, by \rf{48}, \rf{49} and a change of variables, one has
${\fk I}^+(\bfc)= a^{-1}{\fk I}^+(\bfc')$. Also, an inspection of \rf{47}, \rf{Ep} and \rf{Cong}
reveals that ${\fk S}(\bfc)= a {\fk S}(\bfc')$ so that
\be{inprim} {\fk I}^+(\bfc) {\fk S}(\bfc)= {\fk I}^+(\bfc') {\fk S}(\bfc'). \ee

\subsection{Positive solutions}
Throughout this section we continue to use notational conventions from earlier sections of this chapter as well as 
those from Chapter 3. The object is to establish the asymptotic relation \rf{32+}. 
We consider $d,k$ as fixed and suppose that $s>n_0(d)$. For
any non-zero integers $c_j$, we write
$$ F(\al ) = f(c_1\al) f(c_2\al) \cdots f(c_s\al) $$
where $f(\al)$ is the exponential sum \rf{321}. Then, by orthogonality, 
$$
\MM^+_\bfc (\bfX) = \int_0^1 F(\al)\,d\al. 
$$

We evaluate $\MM^+_\bfc (\bfX)$ by the 
Hardy-Littlewood method. In doing so, we may suppose that \rf{322}
holds (by symmetry), Also, we suppose now that $(c_1;\ldots;c_s)=1$. 


Let $P$ be as in \rf{defP}, and let $U$ be the parameter introduced in Lemma \ref{weyl2}. Now put
$$ W = U^{s} . $$
 Let $\fk K$
denote the disjoint union of the intervals
$$ \{ \al\in[0,1]: |\al -a/q| \le WP^{-d}\} $$
with $0\le a\le q\le W$ and $(a,q)=1$. Let ${\fk k} = [0, 1]\setminus \fk K$. 
For $\fk a$ one of the sets $\fk K$, $\fk k$, we write
$$ I({\fk a}) = \int_{\fk a} F(\al)\,d\al  $$
and then have
\be{353} \MM^+_\bfc (\bfX) = I({\fk K}) + I({\fk k}).\ee

For $1\le j\le s$, consider the sets
$$ {\fk l}_j =\{ \al\in[0,1]: \, c_j\al \in {\fk n}(U)\}. $$
Suppose that $\al\in[0,1]$ is in none of the ${\fk l}_j$ $(1\le j\le s)$.
Then, for all $j\le s$, one has $c_j\al\in{\fk N}(U)$, and this shows that there
are coprime $b_j\in\ZZ$, $r_j\in\NN$ with $|c_j\al- (b_j/r_j)| \le r_j^{-1}
UP^{-d}$ and $r_j\le U$. We now compare the various approximations
$$ \Big| \al - \f{b_j}{c_j r_j}\Big| \le \f{U}{|c_j|r_j P^d}.$$
For $1\le i<j\le s$ one has
$$  \Big|    \f{b_j}{c_j r_j}-\f{b_i}{c_i r_i} \Big| \le \f{U}{P^d}\Big( \f1{|c_j|r_j} + \f1{|c_i|r_i}\Big), $$
whence
$$ 
| c_ir_i b_j - c_jr_j b_i | \le UP^{-d} (|c_i|r_i + |c_j|r_j) \le 2U^2P^{-d}
|\bfc|< 1
$$ 
provided only that $|\bfc|\le P^{1/2}$, as we temporarily assume. Consequently,
we may write $b_j/(c_j r_j) = a/q$ with $(a,q)=1$, for all $1\le j\le s$.
Let $p^\nu\parallel q$, with $\nu\ge 1$. Then, since $(c_1;\ldots;c_s)=1$,
there is a $j$ with $p\nmid c_j$. But $q|c_j r_j$ then implies $p^\nu| r_j$.
It follows that $q|r_1r_2\cdots r_s$, and hence that $q\le U^s$. This
shows that $\al\in \fk K$, and we may conclude that
$$ I({\fk k}) \le \sum_{j=1}^s  \int_{{\fk l}_j} |F(\al)|\,d\al. $$
An inspection of Lemma \ref{weyl2} shows that there is a $\delta>0$
depending only on $d$, $k$ and $s$, and such that 
$$ \sup_{\al\in{\fk l}_1} |f(c_1\al)| \le \sup_{\gamma\in{\fk n}(U)}
|f(\gamma)| \ll PX_k^{-2\delta}. $$
Hence, on writing $\sgm=s-\f12$, H\"older's inequality yields
$$ \int_{{\fk l}_1} |F(\al)|\,d\al \ll P^{1/2}X_k^{-\delta}
\Big(\int_0^1 |f(c_1\al)|^\sgm\,d\al\Big)^{1/(2\sgm)}
\prod_{j=2}^s \Big(\int_0^1 |f(c_j\al)|^\sgm\,d\al\Big)^{1/\sgm}. $$
An obvious substitution and Lemma  \ref{L341} yield
$$ \int_0^1 |f(c_j\al)|^\sgm\,d\al = 
\int_0^1 |f(\al)|^\sgm\,d\al \ll P^{\sgm-d},$$
which combines with the previous inequality to 
$$ \int_{{\fk l}_1} |F(\al)|\,d\al \ll P^{s-d} X_k^{-\delta}. $$
By symmetry, the same bound holds with ${\fk l}_1$ replaced by any other 
${\fk l}_j$, so that may now conclude that
\be{minor}
I({\fk k}) \ll P^{s-d} X_k^{-\delta}. \ee

We may now concentrate on the major arcs $\fk K$. For $\al\in\fk K$, there is
a unique pair $a,q$ of coprime integers with $0\le a \le q$, $1\le q\le W$
and $\al=(a/q)+\beta$ satisfying $|\beta|\le WP^{-d}$. Lemma \ref{L322}
gives
\be{errf} f(c_j\al) = q^{-k}S(q,ac_j)v(c_j\beta) + O(E(q,\bfX) (|c_j|W)^k). \ee
For the leading term on the right hand side, crude use of Lemmas \ref{L323} and \ref{L324}
supply the trivial upper bound $q^{-k}S(q,ac_j)v(c_j\beta)\ll P$, and we also have
$f(c_j\al)\ll P$. Hence, for the difference in \rf{errf}, we have the alternative
yet trivial bound $O(P)$. We may now multiply 
together
to infer that
\be{47N} F(\al) = q^{-ks} \prod_{j=1}^s S(q,ac_j) v(c_j\beta) + O(P^{s-1} E(q,\bfX) |\bfc|^k W^k).
\ee
By \rf{323}, whenever $q\le W$ one has $E(q,\bfX) \ll PWX_k^{-1}$, and the measure of $\fk K$ is $O(W^3 P^{-d})$. Hence, on integrating \rf{47N} over $\fk K$, one infers that
$$ I({\fk K}) = {\fk S}_k(\bfc, W) \int_{-WP^{-d}}^{WP^{-d}} v(c_1\beta)\ldots v(c_s\beta) \,d\beta + O(P^{s-d}W^{k+4}X_k^{-1} |\bfc|^k ) $$
where 
${\fk S}_k(\bfc, W)$ is given by \rf{45}.
Within the integral on the right hand side, we use \rf{3210}, then substitute $\beta$ for $P^k\beta$ and recall \rf{411} to recast the previous display in the form
\be{48N} I({\fk K}) = {\fk S}_k(\bfc, W) {\fk I}_k^+(\bfc, W) P^{s-d}+ O(P^{s-d}W^{k+4}X_k^{-1} |\bfc|^k ). \ee
By \rf{44}, \rf{46}, \rf{410} and \rf{412}, 
$$ {\fk S}_k(\bfc, W) {\fk I}_k^+(\bfc, W) = {\fk S}_k(\bfc) {\fk I}_k^+(\bfc)
+ O( W^{-1/(2d)} |\bfc|^{s/d}). $$
We now define
\be{C+} {\fk E}^+(\bfc) =  {\fk S}_k(\bfc) {\fk I}_k^+(\bfc), \ee
and then first deduce from \rf{48N} that
$$  I({\fk K}) ={\fk E}^+(\bfc)P^{s-d} +O(P^{s-d} X_k^{-\delta} |\bfc|^{s+k}) $$
holds for some sufficiently small positive $\delta$, and further, from \rf{minor} and \rf{353}
that indeed \rf{32+} is valid. Along the way we have assumed that $|\bfc|\le P^{1/2}$.
However, in the contrary case, the condition that $s>2d$ implies that
 $P^{s-d} X_k^{-\delta} |\bfc|^{s+k} \gg P^{s}$, which shows that again \rf{32+} holds, this time for trivial reasons. This completes the proof of \rf{32+} when $\bfc$ is primitive.

Now suppose that $\bfc$ is not primitive, and that $\bfc= a \bfc'$ with $a\in\NN$ and $\bfc'$ primitive.
Then it is immediate that $ \MM^+_\bfc (\bfX)=\MM^+_{\bfc'} (\bfX) $, and we may apply \rf{32+} to
$\MM^+_{\bfc'} (\bfX)$ to derive \rf{32+} for $\MM^+_\bfc (\bfX)$ with ${\fk E}^+(\bfc)={\fk E}^+(\bfc')$.
This establishes \rf{32+} in all cases. In passing we mention that \rf{C+} holds also in
the case where $\bfc$ is not primitive, as one finds from \rf{inprim}.

\smallskip

All other conclusions concerning $\MM^+_\bfc (\bfX)$ in Theorem \ref{S4.1} are also available.
By Lemmas \ref{Lss} and \ref{Lsi} we see that ${\fk E}^+(\bfc)$ is real and non-negative. 
Moreover, whenever \rf{33} has non-trivial solutions in $\QQ_p$ for all primes $p$, and a
solution in positive real numbers, then the $c_j$ cannot all be of the same sign, and  Lemmas \ref{Lss} and \ref{Lsi} show ${\fk E}^+(\bfc)>0$. Finally, if \rf{33} fails to have non-trivial solutions in some
$\QQ_p$, or in positive reals, then in particular there are no solutions of \rf{33} in natural numbers.
Hence $\MM^+_\bfc (\bfX)=0$, which is compatible with \rf{32+} only when ${\fk E}^+(\bfc)=0$. This completes the
proof of Theorem \ref{S4.1} for $\MM^+_\bfc (\bfX)$.

\subsection{The proof of Theorem 4.1 completed}
It remains to establish Theorem \ref{S4.1} for $\MM_\bfc (\bfX)$. When $d$ is even, an inspection of the
definitions of $\MM_\bfc (\bfX)$ and $\MM^+_\bfc (\bfX)$ reveals that $\MM_\bfc (\bfX)=2^{ks}\MM^+_\bfc (\bfX)$, and
all conclusions concerning $\MM_\bfc (\bfX)$ follow from those for $\MM^+_\bfc (\bfX)$ if we put
${\fk E}(\bfc)=2^{ks}{\fk E}^+(\bfc)$.

When $d$ is odd, a similar but more elaborate argument applies. Since $-1$ is a $d$-th power, we have
\be{430}\MM_\bfc (\bfX) 
= 2^k \multsum{\eta_j\in\{1,-1\}}{1\le j\le s} \MM^+_{\bm\eta\bfc}(X)
\ee
where $\bm\eta\bfc = (\eta_1 c_1,\ldots,\eta_s c_s)$ is the coordinate product. Since $-1$ is also a $d$-th power residue, modulo $q$, it follows from \rf{324} that $S(q,-a)=S(q,a)$, and consequently, via \rf{singser},
one readily confirms that ${\fk S}(\bm\eta\bfc)={\fk S}(\bfc)$ for all $\bm\eta$ that occur in the sum \rf{430}. With the asymptotic formula \rf{32+} already in hand, the cognate expansion \rf{32} with
$$ {\fk E}(\bfc) = 2^k {\fk S}(\bfc) \multsum{\eta_j\in\{1,-1\}}{1\le j\le s} {\fk I}^+(\bm\eta\bfc) 
$$
is immediate consequence of \rf{430}. We also see that ${\fk E}(\bfc)$ is real and non-negative. We already showed that ${\fk S}(\bfc)>0$ holds whenever \rf{33} has non-trivial $p$-adic solutions for all primes $p$.
Further, since $d$ is now odd, the equation \rf{33} has a real solution with all $y_j$ non-zero. For this
solution,
put $\eta_j=y_j/|y_j|$. The analogue of \rf{33} with $\bm\eta\bfc$ in place of $\bfc$ has the solution
$(|y_1|,\ldots,|y_s|)$, whence Lemma \ref{Lsi} yields ${\fk I}^+(\bfc)>0$. We conclude that ${\fk E}(\bfc)>0$ holds whenever \rf{33} has non-trivial solutions in all $\QQ_p$. If that fails, then \rf{32} implies     
 ${\fk E}(\bfc)=0$. The proof of Theorem \ref{S4.1} is now complete.

We remark that  ${\fk E}(\bfc)$
can be written in a form that is independent of the parity of $d$. In fact, if $d$ is odd, one deduces from \rf{329} that
$$ 2^k \multsum{\eta_j\in\{1,-1\}}{1\le j\le s} {\fk I}^+(\bm\eta\bfc) = {\fk I}(\bfc)
$$
where
\be{432}
{\fk I}(\bfc) = \int_{-\infty}^{\infty} \int_{[0,1]^{ks}} e(\beta(c_1\la\bft_1\ra^d+\ldots+c_s\la\bft_s\ra^d))\,\dd\bft\,\dd\beta.
\ee
Note that ${\fk I}(\bfc)$ is also defined when $d$ is even, and in this case, an inspection of
\rf{329} and \rf{432} show that $2^{ks} {\fk I}^+(\bfc)={\fk I}(\bfc)$. It follows that \rf{32}
holds for all values of $d$ with 
\be{433} {\fk E}(\bfc)={\fk S}(\bfc){\fk I}(\bfc). \ee

\subsection{Primitive solutions} In this section, we study the equation
\rf{31} in ``transposed'' form. Note that $\MM_\bfc(\bfX)$ equals the number of solutions of
\be{440}  \sum_{j=1}^s c_j (x_{1,j}x_{2,j} \cdots x_{k,j})^d = 0
\ee
with $1\le |x_{i,j}|\le X_i$ for $1\le i\le k$, $1\le j\le s$. We require an asymptotic formula for the number $\MM^*_\bfc(\bfX)$ of those solutions counted by $\MM_\bfc(\bfX)$ that satisfy the
additional constraints
$$  (x_{i,1};x_{i,2}; \ldots; x_{i,s})=1 \quad (1\le i\le k). $$
With this in view, one may arrange the solutions of \rf{440} according to the
values of $l_i= (x_{i,1};x_{i,2}; \ldots; x_{i,s})$ to infer that
$$ \MM_\bfc(\bfX) = \sum_{\bfl\le \bfX} \MM^*_\bfc(X_1/l_1,\ldots,X_k/l_k). $$
By one of M\"obius' inversion formulae, we deduce that
$$ \MM^*_\bfc(\bfX) = \sum_{\bfl\le \bfX}\mu(\bfl) \MM_\bfc(X_1/l_1,\ldots,X_k/l_k) $$
where in the interest of brevity we put $\mu(\bfl)=\mu(l_1)\mu(l_2)\cdots\mu(l_k)$.
We now suppose that the hypotheses in Theorem \ref{S4.1} hold, and inject the asymptotic
formula \rf{32} into the preceding identity. This yields
$$ \MM^*_\bfc(\bfX) = {\fk E}(\bfc) \la\bfX\ra^{s-d}\sum_{\bfl\le \bfX}\frac{\mu(\bfl)}{\la\bfl\ra^{s-d}}
+ O ( \la\bfX\ra^{s-d} (\min X_i)^{-\delta/2} |\bfc|^{s+k} ),   $$
as one readily confirms. Routine estimates also show 
$$ \sum_{\bfl\le \bfX}\frac{\mu(\bfl)}{\la\bfl\ra^{s-d}} = \zeta(s-d)^{-k} +O((\min X_i)^{-1}),$$
so that we may  conclude as follows.

\begin{lem}\label{primsol} Let $s>n_0(d)$. Then there is a positive number $\delta'$
such that whenever $c_j\in\ZZ\setminus\{0\}$ one has
$$  \MM^*_\bfc(\bfX) = \zeta(s-d)^{-k}{\fk E}(\bfc) \la\bfX\ra^{s-d} +O(\la\bfX\ra^{s-d} (\min X_i)^{-\delta'} |\bfc|^{s+k} ).   $$
\end{lem}

\section{Synthesis}
The scene is prepared for a swift derivation of Theorem \ref{Thm01}. The starting point is \rf{010}, and the strategy is to show that the function $\theta(\bfm)$ is part of a suitable family for Theorem \ref{S2.1} to deliver Theorem \ref{Thm01}. From now on, we consider $\bfa=(a_0,\ldots,a_s)$
as fixed, once and for all, and suppose that the hypotheses of Theorem \ref{Thm01} are satisfied.

We begin by observing that the equations \rf{01} and \rf{440} become identical if one takes $s=n+1$
and $c_j=a_{j-1}$. An examination of the definition of $\theta(\bfm)$ and \rf{011} now reveals that
$\Theta(\bfX)=\MM^*_\bfa(\bfX)$. Hence, by Lemma \ref{primsol}, there is a positive $\delta$ such that
$$ \Theta(\bfX) =  \zeta(n+1-d)^{-k}{\fk E}(\bfa) \la\bfX\ra^{n+1-d}
+ O(\la\bfX\ra^{n+1-d} (\min X_i)^{-\delta}) $$
holds, and we conclude that $\theta$ obeys condition (I) with
\be{defC} c = \zeta(n+1-d)^{-k}{\fk E}(\bfa). \ee

Next, we check condition (II). 
Fix $r$ with $1\le r\le k-1$ and put $l=k-r$. With $\bfu\in\NN^r$ and $\bfV\in [1,\infty)^l$, we have to evaluate the sum
\be{51} \Theta_\bfu(\bfV) = \sum_{\bfv\le \bfV} \theta (\bfu,\bfv). \ee
Progress depends on a diophantine interpretation of this sum that we prepare by rewriting the equation \rf{01} in a notation more suitable for the current needs. Thus we consider
$$ \sum_{j=0}^n  c_j(y_{1,j}\ldots y_{r,j} z_{1,j}\ldots z_{l,j})^d =0 $$
 and observe that $\theta(\bfu,\bfv)$ is the number of its solutions in primitive vectors $\bfy_i, \bfz_{i'}\in\ZZ^{n+1}$
with non-zero coordinates and $|\bfy_i|=u_i$, $|\bfz_{i'}|=v_{i'}$ $(1\le i\le r,\, 1\le i'\le l)$. For a fixed permissible choice of
$\bfy_1,\ldots,\bfy_r$ we may sum over $\bfv\le \bfV$. One then obtains a quantity examined in Lemma \ref{primsol}, but
with $l$ in place of $k$, with $s=n+1$ and $\bfc=\bfc(\bfy_1,\ldots,\bfy_r)$ given by
$$c_{j+1} = a_j( y_{1,j}\ldots y_{r,j} )^d\quad  (0\le j\le n). $$
Consequently, on writing
$$ {\cal Y} (\bfu) =\{(\bfy_1,\ldots,\bfy_r): \bfy_i\in\ZZ^{n+1} \text{ primitive}, |\bfy_i|=u_i\, (1\le i\le r)\}, $$
we find that
$$ \Theta_\bfu(\bfV) = \sum_{(\bfy_1,\ldots,\bfy_r)\in  {\cal Y} (\bfu)} \MM^*_{\bfc(\bfy_1,\ldots,\bfy_r)} (\bfV). $$
By Lemma \ref{primsol}, we now infer that
\be{53} 
 \Theta_\bfu(\bfV) = \zeta(n+1-d)^{-k} \la\bfV\ra^{n+1-d} \sum_{(\bfy_1,\ldots,\bfy_r)\in  {\cal Y} (\bfu)}
{\fk E} (\bfc(\bfy_1,\ldots,\bfy_r)) + E
\ee
where
$$ E \ll \la\bfV\ra^{n+1-d} (\min V_i)^{-\delta}  \sum_{(\bfy_1,\ldots,\bfy_r)\in  {\cal Y} (\bfu)}
 |\bfc(\bfy_1,\ldots,\bfy_r)|^{n+1+k}. $$
Note that $| y_{1,j}\ldots y_{r,j}| \le |\bfy_1|\cdots |\bfy_r| = \la\bfu\ra$. This yields 
$$ |\bfc(\bfy_1,\ldots,\bfy_r)| \ll \la\bfu\ra^d. $$
Further, there are now more that $O(\la\bfu\ra^{nr})$ elements in $ {\cal Y} (\bfu)$. It follows that
\be{54} E\ll \la\bfV\ra^{n+1-d} (\min V_i)^{-\delta}  \la\bfu\ra^{nr+d(n+1+k)} \ll
\la\bfV\ra^{n+1-d} (\min V_i)^{-\delta} |\bfu|^D \ee
where $D= r(nr+d(n+1+k))$. By \rf{51}, \rf{53} and \rf{54}, we  see that $\theta$ satisfies condition (II).

For condition (III), we note that $\theta_\sigma=\theta$ holds for all $\sigma\in S_k$ by symmetry, as one
confirms from \rf{01}. We have now proved that the 
function $\theta$ alone is a $(n+1-d,c,D,1,\delta)$-family where $c$ and $D$ are as above, and $\delta$ is a sufficiently small positive number. We may now apply Theorem \ref{S2.1} with $N=B^{1/(n+1-d)}$.
Then,  whenever $c>0$, the conclusions of Theorem \ref{Thm01} follow from \rf{010},
and by  \rf{defC} we also see that the constant $C$ is given by
\be{final}  C= \frac{{\fk S}(\bfa) {\fk I}(\bfa)}{2^k(k-1)!\,\zeta(n+1-d)^k}. \ee
 If $c=0$, then by \rf{final}, we have $C=0$, and   by \rf{defC}, \rf{433} 
and Lemmas \ref{Lss} and \ref{Lsi}, the equation \rf{04} has only the trivial solution in at least one of
the fields $\QQ_p$ or $\RR$. In this case $\mathrm N (B)=0$, and the conclusions of Theorem \ref{Thm01} again follow, with the same formula \rf{final} for $C$. 
This completes the proof of Theorem \ref{Thm01}, and from \rf{final}, \rf{47} and  \rf{Ep} we infer that $C$ is the product of local densities that the Hardy-Littlewood method predicted.

\end{document}